\documentclass[12pt,leqno]{article}
\usepackage{amssymb}
\usepackage[mathscr]{eucal}
\usepackage{amsmath,amssymb,latexsym,theorem,bbm}
\usepackage{color}
\newfont{\msa}{msam10 scaled\magstep1}
\newfont{\ssmsa}{msam9}
\newfont{\smsa}{msam10}
\newfont{\sms}{msbm10}
\newfont{\sseufb}{eufb9}
\newfont{\seufb}{eufb10}
\newfont{\eufb}{eufb10 scaled\magstep1}
\newfont{\eusb}{eusb10 scaled\magstep1}
\newfont{\hcmr}{cmr17 scaled\magstep5}

\newcommand{\DS}{\displaystyle}

\newcommand{\EE}{\mathsf{E}}

\newcommand{\NN}{\mathbb{N}}
\newcommand{\PP}{\mathsf{P}}

\newcommand{\RR}{\mathbb{R}}

\newcommand{\sSS}{\raise-0.5truemm\hbox{\sms S}}

\newcommand{\cA}{{\mathcal A}}

\newcommand{\cF}{{\cal F}}

\newcommand{\dd}{\mathsf{d}}

\newcommand{\ee}{\mathsf{e}}
\newcommand{\var}{\mathsf{Var\,}}

\renewcommand{\leq}{\leqslant}
\renewcommand{\geq}{\geqslant}

\newcommand{\proofend}{\hfill\mbox{$\Box$}}

\numberwithin{equation}{section}

\theoremstyle{change} \theorembodyfont{\em}
\newtheorem{Lem}{Lemma.}[section]
\newtheorem{Thm}[Lem]{Theorem.}
\newtheorem{Pro}[Lem]{Proposition.}
\newtheorem{Cor}[Lem]{Corollary.}
\newtheorem{Def}[Lem]{Definition.}
\theorembodyfont{\rm}
\newtheorem{Rem}[Lem]{Remark.}
\newtheorem{Ex}[Lem]{Example.}

\begin{document}

\title{\vspace*{-15mm}
       \bfseries\Large General $\alpha$-Wiener bridges}

\author{{\sc M\'aty\'as Barczy} and {\sc Peter Kern}}

\date{}

\maketitle



\renewcommand{\thefootnote}{}
\footnote{M. Barczy (corresponding author),
             Faculty of Informatics,
             University of Debrecen,
             Pf.12,
             H--4010 Debrecen,
             Hungary,
             Email: barczy.matyas@inf.unideb.hu}
\footnote{P. Kern,
             Mathematical Institute,
             Heinrich-Heine-Universit$\rm{\ddot{a}}$t D\"usseldorf,
             Universit$\rm{\ddot{a}}$tsstr.~1,
             D-40225 D\"usseldorf,
             Germany,
             Email: kern@math.uni-duesseldorf.de}
\footnote{The first author has been supported by the Hungarian
  Scientific Research Fund under Grant No.\ OTKA T-079128.
This work has been finished while M. Barczy was on a post-doctoral position
 at the Laboratoire de Probabilit\'es et Mod\`{e}les Al\'eatoires,
 University Pierre-et-Marie Curie, thanks to NKTH-OTKA-EU FP7 (Marie Curie action)
 co-funded 'MOBILITY' Grant No. OMFB-00610/2010. }


\begin{abstract}
An \ $\alpha$-Wiener bridge is a one-parameter generalization of the usual Wiener bridge, where the parameter \ $\alpha>0$ \ represents a mean reversion force to zero.
We generalize the notion of \ $\alpha$-Wiener bridges
 to continuous functions \ $\alpha:[0,T)\to\RR.$
We show that if the limit \ $\lim_{t\uparrow T}\alpha(t)$ \ exists and is positive, then a general
 \ $\alpha$-Wiener bridge is in fact a bridge in the sense that it converges to \ $0$ \ at time \ $T$ \ with
 probability one.
Further, under the condition \ $\lim_{t\uparrow T}\alpha(t)\ne1$ \
we show that the law of the general \ $\alpha$-Wiener bridge can not
coincide with the law of any non time-homogeneous Ornstein-Uhlenbeck type bridge.
In case \ $\lim_{t\uparrow T}\alpha(t)=1$ \ we determine all the
 Ornstein-Uhlenbeck type processes from which one can derive the general \ $\alpha$-Wiener bridge by
 conditioning the original Ornstein-Uhlenbeck type process to be in zero at time \ $T$.
\end{abstract}

\section{Introduction}

This paper deals with the so-called $\alpha$-Wiener bridges.
Let \ $T\in(0,\infty)$ \ be fixed.
For all \ $\alpha\in\RR$, \ let us consider the stochastic differential equation (SDE)
 \begin{align}\label{alpha_W_bridge}
  \begin{cases}
   \dd X_t=-\frac{\alpha}{T-t}\,X_t\,\dd t+\dd B_t,\qquad t\in[0,T),\\
   \phantom{\dd} X_0=0,
  \end{cases}
 \end{align}
 where \ $(B_t)_{t\geq 0}$ \ is a one-dimensional standard Wiener process defined
 on a filtered probability space \ $(\Omega,\cA,(\cA_t)_{t\geq 0},\PP)$ \
 satisfying the usual conditions (the filtration being constructed by the help of \ $B$),
 i.e., \ $(\Omega,\cA,\PP)$ \ is complete, \ $(\cA_t)_{t\geq 0}$ \ is right continuous,
 \ $\cA_0$ \ contains all the $\PP$-null sets in \ $\cA$ \ and \ $\cA_\infty=\cA$, \ where
 \ $\cA_\infty:=\sigma\left(\bigcup_{t\geq 0}\cA_t\right)$, \ see, e.g., Karatzas and Shreve
 \cite[Section 5.2.A]{KarShr}.
By {\O}ksendal \cite[Theorem 5.2.1]{Oks} or Jacod and Shiryaev \cite[Chapter III, Theorem 2.32]{JacShi},
 the SDE \eqref{alpha_W_bridge} has a unique strong solution, namely,
 \begin{align}\label{alpha_W_bridge_intrep}
   X_t=\int_{0}^t\left(\frac{T-t}{T-s}\right)^\alpha\,\dd B_s,\qquad t\in[0,T),
 \end{align}
 as it can be checked by It\^{o}'s formula.
The Gauss process \ $(X_t)_{t\in[0,T)}$ \ given by \eqref{alpha_W_bridge_intrep} is called an
 $\alpha$-Wiener bridge.
 More generally, we call any almost surely continuous (Gauss) process on the time interval
 \ $[0,T)$ \ having the same finite-dimensional distributions as \ $(X_t)_{t\in[0,T)}$ \
 an \ $\alpha$-Wiener bridge.
To our knowledge, these kinds of processes have been first considered in the case of \ $\alpha>0$
 \ by Brennan and Schwartz \cite{BreSch}; see also Mansuy \cite{Man}.
In Brennan and Schwartz \cite{BreSch}, \ $\alpha$-Wiener bridges (with \ $\alpha>0$)
 \ are used to model the arbitrage profit associated with a given futures contract in the absence
 of transaction costs.
Sondermann, Trede and Wilfling  \cite{SonTreWil} and Trede and Wilfling \cite{TreWil}
use the SDE \eqref{alpha_W_bridge} (with \ $\alpha>0$) to describe the fundamental component
 of an exchange rate process and they call the process $X$ a scaled Brownian bridge.
The essence of these models is that
 the coefficient of $X_t$ in the drift term in \eqref{alpha_W_bridge} represents some kind of mean
 reversion, a stabilizing force that keeps pulling the process towards its mean (i.e., to zero),
 and the absolute value of this force is  increasing proportionally to the inverse of the remaining time $T-t,$
 with the rate constant $\alpha$.
 Note also that in case of \ $\alpha=1$ \ the process \ $(X_t)_{t\in[0,T)}$ \ is nothing else
 but the usual Wiener bridge (from $0$ to $0$ over \ $[0,T]$).

It is known that in case of \ $\alpha>0$, \ the $\alpha$-Wiener bridge \ $(X_t)_{t\in[0,T)}$ \ given by
 \eqref{alpha_W_bridge_intrep} has an almost surely continuous extension \ $(X_t)_{t\in[0,T]}$ \ to the
 time-interval \ $[0,T]$ \ such that \ $X_T=0$ \ with probability one, see, e.g.,
 Mansuy \cite[page 1023]{Man} or Barczy and Pap \cite[Lemma 3.1]{BarPap2}.
For positive values of \ $\alpha$, \ the possibility of such an extension is based on a strong law of
 large numbers for square integrable local martingales.
In case of \ $\alpha\leq 0$, \ there does not exist an almost surely continuous extension of the process
 \ $(X_t)_{t\in[0,T)}$ \ onto \ $[0,T]$ \ which would take some constant at time \ $T$ \ with probability one
 (i.e., which would be a bridge).
Indeed, for \ $\alpha=0$ \ the process \ $X$ \ is nothing else but a standard
 Wiener process (which is not a constant at time \ $T$ \ with probability one),
 and in case of \ $\alpha<0$ \ it can be checked that the second moment of \ $X_t$ \
 (given by \eqref{alpha_W_bridge_intrep}) converges to infinity as \ $t\uparrow T$.
Hence in case of \ $\alpha<0$ the assumption of the existence of an almost surely continuous
 extension to \ $[0,T]$ \ such that this extension takes some constant at time \ $T$ \ with probability
 one (i.e., we have a bridge) would result in a contradiction.
We note that another proof of the impossibility of such an extension in the case of \ $\alpha<0$ \
 can be found in Barczy and Pap \cite[Remark 3.5]{BarPap2}.
For a detailed discussion of sample path properties of $\alpha$-Wiener bridges, see Barczy and Pap \cite{BarPap2}.
 Finally we remark that an $\alpha$-Wiener bridge (for all \ $\alpha\in\RR$) \ can be represented
 as a space-time transformed Wiener process, for a detailed discussion
 see Barczy and Igl\'oi \cite[Remark 2.4]{BarIgl}.

The main contribution of the present paper is a detailed study of the question of
 so-called identical bridges for \ $\alpha$-Wiener bridges.
Up to our knowledge these kinds of investigations were started by Benjamini and Lee \cite{BenLee}
 for usual Wiener bridges.
 They determined all the one-dimensional diffusion processes \ $(Y_t)_{t\geq 0}$ \ being weak
 solutions of the SDE
\begin{align}\label{SEGED2}
  \dd Y_t=b(Y_t)\,\dd t+\,\dd B_t,\qquad t\geq 0,
\end{align}
 which satisfy the following property: for any \ $T>0$ \ and any \ $x\in\RR$, \
  the bridge from $x$ to $x$ over $[0,T]$ derived from $Y$
 is the Wiener bridge from \ $x$ \ to \ $x$ \ over \ $[0,T]$.
Namely, under the condition that the function \ $b:\RR\to\RR$ \ is bounded and twice continuously
 differentiable they showed that either \ $b$ \ is constant or
 $$
   b(x)=k\tanh(kx+c),\qquad x\in\RR,\quad k,c\in\RR.
 $$
This result has been extended by Fitzsimmons \cite{Fit} in two ways.
Firstly, he studied bridges constructed from more general time-homogeneous Markov processes
 with values in an abstract state space under suitable regularity conditions.
Secondly, under some additional continuity condition, he showed that
 if \ $X$ \ and \ $Y$ \ are time-homogeneous Markov processes (with values in an abstract state space)
 and there exist real numbers \ $x_0,y_0\in\RR$ \ and \ $T_0>0$ \ such that the law of the bridge
 from \ $x_0$ \ to \ $y_0$ \ over \ $[0,T_0]$ \ derived from \ $X$ \ coincides with the law of the
 bridge from \ $x_0$ \ to \ $y_0$ \ over \ $[0,T_0]$ \ derived from \ $Y$, \ then the same statement
 holds for bridges from \ $x$ \ to \ $y$ \ over \ $[0,T]$ \ with arbitrary \ $x,y\in\RR$ \ and
 \ $T>0$.
\ Recently, Borodin \cite{Bor} considered the original question of Benjamini and Lee \cite{BenLee}
 replacing the usual Wiener bridge by the Bessel bridge or the (radial) Ornstein-Uhlenbeck bridge.

In Section \ref{Section_alpha_Wiener} we generalize the notion of \ $\alpha$-Wiener bridges.
Namely, for \ $T>0$ \ and a continuous function \ $\alpha:[0,T)\to\RR$ \ we consider the SDE
  \begin{align}\label{gen_alpha_W_bridge}
  \begin{cases}
   \dd X_t=-\frac{\alpha(t)}{T-t}\,X_{t}\,\dd t+\dd B_t,\qquad t\in[0,T),\\
   \phantom{\dd} X_0=0,
  \end{cases}
 \end{align}
 where \ $(B_t)_{t\geq 0}$ \ is a one-dimensional standard Wiener process.
This SDE has a unique strong solution \ $(X_t)_{t\in[0,T)}$ \ given in
 \eqref{gen_alpha_W_bridge_intrep} which will be called a Wiener bridge
 with continuously varying parameter \ $\alpha$ \ or a general $\alpha$-Wiener bridge.
 More generally, we call any almost surely continuous (Gauss) process on the time interval
 \ $[0,T)$ \ having the same finite-dimensional distributions as \ $(X_t)_{t\in[0,T)}$ \
 a general $\alpha$-Wiener bridge.
In Theorem \ref{THEOREM8} we prove that under the condition that the limit
 \ $\lim_{t\uparrow T}\alpha(t)$ \ exists and is positive,
 we have \ $X_t\to 0$ \ almost surely as \ $t\uparrow T,$ \ which explains that why we can
 use the expression ''bridge'' for \ $X$ \ (at least under the above
 assumption on \ $\alpha$). \ We also examine the question of identical bridges, i.e.,
 whether it is possible to interpret this process as a bridge derived from
 an Ornstein-Uhlenbeck type process  \ $(Z_t)_{t\geq 0}$ \ given by the SDE
 \begin{align}\label{gen_OU_egyenlet}
     \dd Z_t=q(t)\,Z_t\,\dd t + \sigma(t)\,\dd B_t,\qquad t\geq0,
 \end{align}
 with an initial condition \ $Z_0$ \ having a Gauss distribution independent of \ $B$,
 where \ $q:[0,\infty)\to\RR$ \ and \ $\sigma:[0,\infty)\to\RR$ \  are continuous functions
 and \ $(B_t)_{t\geq 0}$ \ is a standard Wiener process.
Here, and in all what follows, by a bridge derived from \ $Z$ \ we mean the
 construction presented in Section \ref{Section_OU_proc_bridge} (summarized in Theorem \ref{THEOREM2}
 and Definition \ref{DEFINITION_bridge}).
Theorems \ref{THEOREM4} and \ref{THEOREM5} give a complete answer to our question in some sense
 (see also Remark \ref{REMARK8_lacks}).
Namely, it turns out that \ $\lim_{t\uparrow T}\alpha(t) = 1$ \ is a necessary condition for
  the existence of such a process \ $Z$ (see Theorem \ref{THEOREM4}),
  but it also turns out that this is not a sufficient condition (see Example \ref{REM_peldak2}).
 In Theorem \ref{THEOREM5}, given a continuously differentiable function \ $\alpha$ \ with
 \ $\lim_{t\uparrow T}\alpha(t) = 1$, \ we determine all the Ornstein-Uhlenbeck type processes
 \ $(Z_t)_{t\geq 0}$ \ (given by the SDE \eqref{gen_OU_egyenlet}) such that (for fixed \ $T>0$)
 \ the law of the bridge from \ $0$ \ to \ $0$ \ over \ $[0,T]$ \ derived from \ $Z$ \ coincides
 with the law of the general \ $\alpha$-Wiener bridge.

 In Section \ref{Section_spec_alpha_Wiener}, besides giving examples and applications of
 our results in Section \ref{Section_alpha_Wiener}, we will further examine in detail the special
 case of $\alpha$-Wiener bridges for constant \ $\alpha\in\RR$. \
Mansuy \cite[Proposition 4]{Man} showed that
 the law of the $\alpha$-Wiener bridge with constant \ $\alpha>0$,
 \ $\alpha\not=1$ \ can not be the same as the law of the bridge derived from a centered Gauss,
 strictly stationary Markov process having almost surely continuous paths.
To complement this result, in Corollary \ref{THEOREM3}
 we show that the law of the \ $\alpha$-Wiener bridge with constant \ $\alpha\in\RR$, $\alpha\not=1$ \
 can not coincide with the law of the bridge derived from any Ornstein-Uhlenbeck type process
 \ $(Z_t)_{t\geq 0}$ \ given by the SDE \eqref{gen_OU_egyenlet}.
Note that the Ornstein-Uhlenbeck type process \ $Z$ \ is a centered Gauss-Markov process,
 but in general not time-homogeneous, for a detailed discussion see also Section \ref{Section_OU_proc_bridge}.
We will also examine what happens in case of \ $\alpha=1.$ \
More precisely, in Theorem \ref{THEOREM6} in the case of \ $\alpha=1$ \  we
 determine all the Ornstein-Uhlenbeck type processes \ $(Z_t)_{t\geq 0}$ \
  (given by the SDE \eqref{gen_OU_egyenlet}) such that for fixed \ $T>0$ \ the law of the bridge
  from \ $0$ \ to \ $0$ \ over \ $[0,T]$ \ derived from \ $Z$ \
 coincides with the law of the $\alpha$-Wiener bridge with \ $\alpha=1$,
 i.e., with the law of the usual Wiener bridge.
We emphasize that the answer to our problem can not be derived from
 Benjamini and Lee \cite{BenLee} or Fitzsimmons \cite{Fit}.
 Indeed, an Ornstein-Uhlenbeck type process given by the SDE
\eqref{gen_OU_egyenlet} is in general not time-homogeneous, while in
\cite{BenLee}, \cite{Bor} or \cite{Fit} the considered processes are
 time-homogeneous Markov processes.

\section{Preliminaries on Ornstein-Uhlenbeck type bridges}\label{Section_OU_proc_bridge}

In this section we recall the notion and properties of Ornstein-Uhlenbeck type bridges to
 such an extent we will need in the following sections. For a more detailed
discussion and for the proofs of the results, see for example Barczy and Kern \cite{BarBec}
 (where one can also find extensions to multidimensional process bridges).

Let us consider the SDE \eqref{gen_OU_egyenlet}.
By Section 5.6 in Karatzas and Shreve \cite{KarShr}, there exists a strong solution of the SDE \eqref{gen_OU_egyenlet},
namely
\begin{align}\label{gen_OU_egyenlet_megoldas}
  Z_t=\ee^{\bar q(t)}\left(Z_0 + \int_0^t \ee^{-\bar q(s)} \sigma(s) \,\dd B_s\right)\quad\text{ with }
       \quad \bar q(t):=\int_{0}^tq(u)\,\dd u,\qquad t\geq 0,
 \end{align}
 and strong uniqueness for the SDE \eqref{gen_OU_egyenlet} holds, see, e.g.,
 Jacod and Shiryaev \cite[Chapter III, Theorems 2.32 and 2.33]{JacShi}.
Here and in what follows in this section we assume that \ $Z_0$ \ has a Gauss distribution
 independent of the Wiener process \ $(B_t)_{t\geq 0}$.
\ Then we may define the filtration \ $(\cA_t)_{t\geq 0}$ \ such that
 \ $\sigma\{Z_0,B_s : 0\leq s\leq t\}\subset\cA_t$ \ for all \ $t\geq 0$, \ see, e.g.,
 Karatzas and Shreve \cite[Section 5.2.A]{KarShr}.

We will call the process \ $(Z_t)_{t\geq 0}$ \ a one-dimensional Ornstein-Uhlenbeck
 process with continuously varying parameters, or a process of Ornstein-Uhlenbeck type.

One can easily derive that for \ $0\leq s<t$ \ we have
 $$
   Z_{t}=\ee^{\bar q(t)-\bar q(s)}Z_{s}+\int_{s}^t\ee^{\bar q(t)-\bar q(u)}\sigma(u)\,\dd B_u.
 $$
Hence, given \ $Z_s=x,$ \ the distribution of \ $Z_t$ \ does not depend on \ $(Z_{r})_{r\in[0,s)}$ \
 which yields that \ $(Z_{t})_{t\geq0}$ \ is a Markov process.
Moreover, for any \ $x\in\RR$ \ and \ $0\leq s<t$ \
 the conditional distribution of \ $Z_t$ \ given \ $Z_s=x$ \ is Gauss with mean
 \ $\ee^{\bar q(t)-\bar q(s)}x$ \ and with variance
 $$
   \gamma(s,t):=\int_{s}^t\ee^{2(\bar q(t)-\bar q(u))}\sigma^2(u)\,\dd u<\infty.
 $$
In what follows we will make the general assumption that
 \begin{align}\label{sigma_assumption}
       \sigma(t)\ne0  \qquad \text{for all \ $t\geq 0$}.
        \end{align}
This guarantees
 that the variance \ $\gamma(s,t)$ \ is positive for all \ $0\leq s<t$.
\ Hence \ $(Z_{t})_{t\geq0}$ \ is a Gauss-Markov process with transition densities
  \begin{equation}\label{gen_OU_densities}
  p_{s,t}^Z(x,y)
  =\frac{1}{\sqrt{2\pi\gamma(s,t)}}
     \exp\left\{-\frac{(y-\ee^{\bar q(t)-\bar q(s)}x)^2}{2\gamma(s,t)}\right\},
       \qquad 0\leq s<t,\quad x,y\in\RR,
 \end{equation}
 and \ $Z$ \ has almost surely continuous paths.

In Barczy and Kern \cite{BarBec}, for fixed \ $T>0$ \ and \ $a,b\in\RR$ \ we constructed
a Markov process \ $(U_t)_{t\in[0,T]}$ \ with initial distribution
\ $\PP(U_0=a)=1$ \ and with transition densities
 \begin{equation}\label{gen_bridge_densities}
    p_{s,t}^U(x,y)=\frac{p_{s,t}^Z(x,y)\,p_{t,T}^Z(y,b)}{p_{s,T}^Z(x,b)},
      \quad x,y\in\RR,\quad 0\leq s<t<T,
 \end{equation}
such that \ $U_{t}\to b=U_T$ \ almost surely and also in \ $L^2$ \ as \ $t\uparrow T$.
The process \ $(U_t)_{t\in[0,T]}$ \ is called a bridge of Ornstein-Uhlenbeck type from
\ $a$ \ to \ $b$ \ over \ $[0,T]$ \ derived from \ $Z$, \ see also Definition \ref{DEFINITION_bridge}.
The construction is based on Theorem 3.1 in Barczy and Kern \cite{BarBec}, which we recall now
 for completeness and for our later purposes.
For the proofs, see Barczy and Kern \cite{BarBec}.

For all $a,b\in\RR$ and $0\leq s\leq t<T$, let us introduce the notations
\begin{equation}\label{gen_OU_exp}
 n_{a,b}(s,t)
  :=\frac{\gamma(s,t)}{\gamma(s,T)}\,\ee^{\bar q(T)-\bar q(t)} b
     +\frac{\gamma(t,T)}{\gamma(s,T)}\ee^{\bar q(t)-\bar q(s)} a,
 \end{equation}
 and
 \begin{equation}\label{gen_OU_var}
  \sigma(s,t):=\frac{\gamma(s,t)\,\gamma(t,T)}{\gamma(s,T)}.
 \end{equation}

\begin{Lem}\label{gen_bridge}
 Let us suppose that condition \eqref{sigma_assumption} holds.
Let $b\in\RR$ and $T>0$ be fixed. Then for all $0\leq s<t<T$ and $x,y\in\RR$ we have
\begin{align*}
 \frac{p_{s,t}^Z(x,y)\,
     p_{t,T}^Z(y,b)}{p_{s,T}^Z(x,b)}
  =\frac{1}{\sqrt{2\pi\sigma(s,t)}}
    \exp\left\{-\frac{\left(y-n_{x,b}(s,t)\right)^2}
                      {2\sigma(s,t)}\right\},
\end{align*}
which is a Gauss density (in $y$) with mean $n_{x,b}(s,t)$ and
 with variance $\sigma(s,t)$.
\end{Lem}

\begin{Thm}\label{THEOREM2}
 Let us suppose that condition \eqref{sigma_assumption} holds.
For fixed $a,b\in\RR$ and $T>0$, let the process $(U_t)_{t\in[0,T)}$ be given by
 \begin{align}\label{gen_OU_integral}
   \begin{split}
    U_t=  n_{a,b}(0,t)
           +\int_{0}^t\frac{\gamma(t,T)}{\gamma(s,T)}\ee^{\bar q(t)-\bar q(s)} \sigma(s)\,\dd B_{s},
          \qquad t\in[0,T).
    \end{split}
 \end{align}
Then for any $t\in[0,T)$ the distribution of $U_t$ is Gauss
 with mean $n_{a,b}(0,t)$ and with variance $\sigma(0,t).$ \ Especially, $U_t\to b$ almost
 surely (and hence in probability) and in $L^2$ as $t\uparrow T$.
Hence the process $(U_t)_{t\in[0,T)}$ can be extended to an almost surely (and hence stochastically) and
 $L^2$-continuous process $(U_t)_{t\in[0,T]}$ with $U_0=a $ and $U_T=b$.
Moreover, $(U_t)_{t\in[0,T]}$ is a Gauss-Markov process
 and for any $x\in\RR$ and $0\leq s<t<T$ the transition density
 $\RR\ni y\mapsto p_{s,t}^U(x,y)$ of $U_t$ given $U_s=x$ is given by
 \begin{align*}
   p_{s,t}^U(x,y)
     = \frac{1}{\sqrt{2\pi\sigma(s,t)}}
    \exp\left\{-\frac{\left(y-n_{x,b}(s,t)\right)^2}
                      {2\sigma(s,t)}\right\},
        \qquad y\in\RR,
 \end{align*}
 which coincides with the density given in Lemma \ref{gen_bridge}.
\end{Thm}

\begin{Def}\label{DEFINITION_bridge}
Let \ $(Z_t)_{t\geq0}$ \ be the linear process given by the SDE \eqref{gen_OU_egyenlet} with an initial
 Gauss random variable $Z_0$ independent of $(B_t)_{t\geq 0}$
 and let us assume that condition \eqref{sigma_assumption} holds.
For fixed \ $a,b\in\RR$ \ and \ $T>0$, \ the process \ $(U_t)_{t\in[0,T]}$ \ defined in Theorem
 \ref{THEOREM2} is called a bridge of Ornstein-Uhlenbeck type from $a$ to $b$ over $[0,T]$ derived
 from \ $Z$.
More generally, we call any almost surely continuous (Gauss) process on the time-interval
 \ $[0,T]$ \ having the same finite-dimensional distributions as \ $(U_t)_{t\in[0,T]}$ \ a bridge
 of Ornstein-Uhlenbeck type from $a$ to $b$ over $[0,T]$ derived from \ $Z$.
\end{Def}

One can also derive a SDE which is satisfied by the Ornstein-Uhlenbeck type bridge,
 see for example Theorem 3.3 in Barczy and Kern \cite{BarBec}.
For completeness and for our later purposes we also recall this result.

\begin{Lem}\label{LEMMA4}
 Let us suppose that condition \eqref{sigma_assumption} holds.
The process \ $(U_t)_{t\in[0,T)}$ \ defined by \eqref{gen_OU_integral} is a unique
 strong solution of the linear SDE
 \begin{align}\label{gen_OU_hid1_egyenlet}
  \begin{split}
   \dd U_t= \left[\left(q(t)-\frac{\ee^{2(\bar q(T)-\bar q(t))}}{\gamma(t,T)}\, \sigma^2(t) \right)\,U_{t}
            + \frac{\ee^{\bar q(T)-\bar q(t)}}{\gamma(t,T)}\, \sigma^2 (t) b\right]\dd t
            + \sigma(t)\,\dd B_t
 \end{split}\end{align}
for \ $t\in[0,T)$ \ and with initial condition \ $U_0=a$, \ and strong uniqueness for the SDE
 \eqref{gen_OU_hid1_egyenlet} holds.
\end{Lem}

Note that an Ornstein-Uhlenbeck type bridge can also be derived using a usual
 conditioning approach, see, e.g., Proposition 3.5 in in Barczy and Kern \cite{BarBec}.
Again, for completeness we also recall this result.

\begin{Pro}
Let \ $a,b\in\RR$ \ and \ $T>0$ \ be fixed. Let \ $(Z_t)_{t\geq 0}$
\ be the linear process given by the SDE \eqref{gen_OU_egyenlet} with initial
 condition \ $Z_0=a$ \ and let us assume that condition \eqref{sigma_assumption} holds.
Let \ $n\in\NN$ \ and \ $0<t_1<t_2<\ldots<t_n<T.$ \ Then the conditional
distribution of \ $(Z_{t_1},\ldots,Z_{t_n})$ \ given \ $Z_T=b$ \
equals the distribution of \ $(U_{t_1},\ldots,U_{t_n}),$ \ where
 \ $(U_t)_{t\in[0,T]}$ \ is an Ornstein-Uhlenbeck type bridge from \ $a$ \
 to \ $b$ \ over \ $[0,T]$ \ derived from \ $Z$.
\end{Pro}

Next we formulate the above presented results in the case of usual Ornstein-Uhlenbeck processes
 and bridges.

\begin{Rem}\label{REMARK3}
In case of \ $q(t)=q\ne0,$ $t\geq 0$, \ and \ $\sigma(t)=\sigma\ne0,$ $t\geq 0$, \
 the bridge of Ornstein-Uhlenbeck type \ $(U_t)_{t\in[0,T]}$ \ from $a$ to $b$ over
 \ $[0,T]$ \ defined in \eqref{gen_OU_integral} has the form
 \begin{align}\label{OU_hid1_intrep}
   U_t=a\,\frac{\sinh(q(T-t))}{\sinh(qT)}
         + b\,\frac{\sinh(qt)}{\sinh(qT)}
         + \sigma  \int_0^t\frac{\sinh(q(T-t))}{\sinh(q(T-s))}\,\dd B_s,
         \qquad t\in[0,T),
 \end{align}
  and admits transition densities
   $$
    p_{s,t}^U(x,y)
       =\frac{1}{\sqrt{2\pi\sigma(s,t)}}
        \exp\left\{-\frac{\left(y-\frac{\sinh(q(t-s))}{\sinh(q(T-s))}\,b
                                 -\frac{\sinh(q(T-t))}{\sinh(q(T-s))}\,x\right)^2}
                    {2\sigma(s,t)}\right\}
   $$
 for all \ $0\leq s<t<T$ \ and \ $x,y\in\RR$, \ where \ $\sigma(s,t)$ \ takes the form
 \[
   \sigma(s,t) = \frac{\sigma^2}{q}\,\frac{\sinh(q(T-t))\sinh(q(t-s))}{\sinh(q(T-s)) }.
 \]
Moreover, the SDE \eqref{gen_OU_hid1_egyenlet} has the form
  \begin{align}\label{OU_hid1_egyenlet}
  \begin{cases}
   \dd U_t=q\left(-\coth(q(T-t))\,U_t+\frac{b}{\sinh(q(T-t))}\right)\,\dd t
               + \sigma \,\dd B_t,\qquad t\in[0,T),\\
   \phantom{\dd} U_0=a,
  \end{cases}
  \end{align}
 and \ $(U_t)_{t\in[0,T)}$ \ given by \eqref{OU_hid1_intrep} is a unique strong solution of this SDE.

In case of \ $q(t)=0$, $t\geq 0$, \ and \ $\sigma(t)=\sigma\ne0,$ $t\geq 0$, \ we get
 \ $\dd Z_t = \sigma\,\dd B_t$, $t\geq 0$, \ and by Section 5.6.B in Karatzas and Shreve \cite{KarShr},
 the Wiener bridge from \ $a$ \ to \ $b$ \ over \ $[0,T]$ \ (derived from \ $Z$)
 given by
 \begin{align}\label{standard_Wiener_bridge}
   \widetilde U_t
      = \begin{cases}
          a + (b-a)\frac{t}{T} + \sigma \int_0^t\frac{T-t}{T-s}\,\dd B_s
              & \text{if \ $t\in[0,T)$,}\\
          b   & \text{if \ $t=T$,}
        \end{cases}
 \end{align}
 is a unique strong solution of the SDE
  \begin{align*}
  \begin{cases}
   \dd \widetilde U_t=\frac{b-\widetilde U_t}{T-t}\,\dd t + \sigma\,\dd B_t,\qquad t\in[0,T),\\
   \phantom{\dd} \widetilde U_0=a.
  \end{cases}
 \end{align*}
\proofend
\end{Rem}

\section{General $\alpha$-Wiener bridges}\label{Section_alpha_Wiener}

In this section first we search for conditions on \ $\alpha$ \ under which the general
 \ $\alpha$-Wiener bridge converges to \ $0$ \ almost surely as \ $t\uparrow T$, \ which will
 explain that why we can use the expression ''bridge'' (at least under the desired conditions on
 \ $\alpha$).
Then we examine whether it is possible to derive general \ $\alpha$-Wiener bridges
 from Ornstein-Uhlenbeck type processes given by the SDE \eqref{gen_OU_egyenlet} by taking a bridge.

Let \ $T>0$ \ be fixed and for a continuous function \ $\alpha:[0,T)\to\RR$ \
 let us consider the SDE \eqref{gen_alpha_W_bridge}.

\begin{Pro}\label{PROPOSITION2}
The SDE \eqref{gen_alpha_W_bridge} has a strong solution given by
\begin{align}\label{gen_alpha_W_bridge_intrep}
   X_t=\int_0^t\exp\left\{-\int_s^t\frac{\alpha(u)}{T-u}\,\dd u\right\}\,\dd B_s,
    \qquad t\in[0,T),
\end{align}
 and strong uniqueness holds for the SDE \eqref{gen_alpha_W_bridge}.
\end{Pro}

\noindent{\bf Proof.}
Since for all \ $S\in[0,T)$, \ the function \ $[0,S]\ni t\mapsto -\frac{\alpha(t)}{T-t}$
 \ satisfies the usual Lipschitz and linear growth conditions, by {\O}ksendal
 \cite[Theorem 5.2.1]{Oks} or Jacod and Shiryaev \cite[Chapter III, Theorem 2.32]{JacShi},
 the linear SDE \eqref{gen_alpha_W_bridge} has a strong solution which is pathwise unique
 (i.e., it has a unique strong solution), and it takes the following form (which can be checked
 by It\^{o}'s formula or see, e.g., Karatzas and Shreve \cite[Section 5.6]{KarShr})
 $$
   X_t=\Phi(t)\int_0^t\Phi^{-1}(s)\,\dd B_s,
    \qquad t\in[0,T),
 $$
 where \ $\Phi(t),$ $t\in[0,T),$ \ is the unique solution of the deterministic
  differential equation (DE)
$$
 \begin{cases}
    \Phi'(t)=-\frac{\alpha(t)}{T-t}\Phi(t),\qquad t\in[0,T),\\
    \Phi(0)=1.
 \end{cases}
$$
 Hence
 $$
   \Phi(t)=\exp\left\{-\int_0^t\frac{\alpha(s)}{T-s}\,\dd s\right\},
     \qquad t\in[0,T),
 $$
and then
\begin{align*}
 X_t&=\exp\left\{-\int_0^t\frac{\alpha(u)}{T-u}\,\dd u\right\}
      \int_0^t\exp\left\{\int_0^s\frac{\alpha(u)}{T-u}\,\dd u\right\}\,\dd B_s\\
    &=\int_0^t\exp\left\{-\int_s^t\frac{\alpha(u)}{T-u}\,\dd u\right\}\,\dd B_s,
    \qquad t\in[0,T).
\end{align*}
\proofend

 We will call the Gauss process \ $(X_t)_{t\in[0,T)}$ \ a Wiener bridge with continuously
 varying parameter \ $\alpha$ \ or a general \ $\alpha$-Wiener bridge.
 More generally, we call any almost surely continuous (Gauss) process on the time interval
 \ $[0,T)$ \ having the same finite-dimensional distributions as \ $(X_t)_{t\in[0,T)}$ \
 a general $\alpha$-Wiener bridge.

\begin{Rem}
 Note that in case of \ $\alpha(t)=\alpha\in\RR,$ $t\in[0,T),$ \ we have
 \begin{align*}
    X_t&=\int_0^t\exp\left\{-\alpha\int_s^t\frac{1}{T-u}\,\dd u\right\}\,\dd B_s
       =\int_0^t\exp\big\{\alpha(\ln(T-t)-\ln(T-s))\big\}\,\dd B_s\\
       &=\int_0^t\left(\frac{T-t}{T-s}\right)^\alpha\,\dd B_s,
       \qquad t\in[0,T),
 \end{align*}
as expected (see \eqref{alpha_W_bridge_intrep}).
 \proofend
\end{Rem}

 In what follows we give a sufficient condition under which the process \ $(X_t)_{t\in[0,T)}$ \
 defined in \eqref{gen_alpha_W_bridge_intrep} has an almost surely continuous extension
 to \ $[0,T]$ \ with \ $X_T=0$, again denoted by  \ $(X_t)_{t\in[0,T]}$. \
 The next theorem is a generalization of Lemma 3.1 in Barczy and Pap \cite{BarPap2} to general
 $\alpha$-Wiener bridges.

\begin{Thm}\label{THEOREM8}
          Let \ $T\in(0,\infty)$ \ be fixed and let \ $(B_t)_{t\geq 0}$ \ be a one-dimensional
          standard Wiener process on a filtered probability space \ $(\Omega,\cF,(\cF_t)_{t\in[0,T)},\PP)$
          \ satisfying the usual conditions, constructed by the help of the standard Wiener process $B$
          (see, e.g., Karatzas and Shreve \cite[Section 5.2.A]{KarShr}).
          If \ $\alpha(T):=\lim_{t\uparrow T}\alpha(t)$ \ exists and \ $\alpha(T)>0$, \ then the process
          \ $(X_t)_{t\in[0,T]}$ \ defined by
   \begin{align}\label{SEGED7}
      \begin{split}
        X_t:=
           \begin{cases}
              \int_0^t\exp\left\{-\int_s^t\frac{\alpha(u)}{T-u}\,\dd u\right\}\,\dd B_s
                 & \text{if \ $t\in[0,T)$,}\\
               0 & \text{if \ $t=T$,}
          \end{cases}
      \end{split}
   \end{align}
 is a centered Gauss process with almost surely continuous paths.
\end{Thm}

\noindent{\bf Proof.}
Due to the fact that the integrand in the stochastic integral of \eqref{SEGED7} is deterministic,
 by Bauer \cite[Lemma 48.2]{Bau}, \ $(X_t)_{t\in[0,T]}$ \ is a centered Gauss process.
To prove almost sure continuity of \ $X$, \ we follow the method of the proof of
Lemma 3.1 in Barczy and Pap \cite{BarPap2}.
For all \ $t\in[0,T)$, \ let
 \[
   M_t:= \int_0^t\exp\left\{\int_0^s\frac{\alpha(u)}{T-u}\,\dd u\right\}\,\dd B_s.
 \]
Then, by Proposition 3.2.10 in Karatzas and Shreve \cite{KarShr},
  \ $(M_t)_{t\in[0,T)}$ \ is a continuous, square-integrable martingale with respect
 to the filtration \ $(\cF_t)_{t\in[0,T)}$ \ and with quadratic variation
 \begin{align*}
   \langle M\rangle_t
     :=\int_0^t\exp\left\{2\int_0^s\frac{\alpha(u)}{T-u}\,\dd u\right\}\,\dd s,\qquad t\in[0,T).
 \end{align*}
If \ $\alpha(T):=\lim_{t\uparrow T}\alpha(t)>0$ \ exists, then for every
 \ $0<\delta_1<\alpha(T)<\delta_2<\delta_1+1/2$ \  one can choose \ $t_0\in(0,T)$ \ such that
\begin{equation}\label{alphbounds}
\delta_1\leq\alpha(t)\leq\delta_2,\qquad t\in[t_0,T].
\end{equation}

 First we consider the case \ $\alpha(T)>1/2$.
 Let \ $\delta_1$ \ and \ $\delta_2$ \ be given such that
 \ $1/2<\delta_1<\alpha(T)<\delta_2<\delta_1+1/2$.
\ Then for all \ $t\in(t_0,T)$ \ we have
 \begin{align}\label{martres}
  \langle M\rangle_t
     & = \int_0^{t_0}\exp\left\{2\int_0^s\frac{\alpha(u)}{T-u}\,\dd u\right\}\,\dd s
         + \int_{t_0}^t\exp\left\{2\int_0^s\frac{\alpha(u)}{T-u}\,\dd u\right\}\,\dd s \nonumber\\
     & = C_1 + \exp\left\{2\int_0^{t_0}\frac{\alpha(u)}{T-u}\,\dd u\right\}
                \int_{t_0}^t\exp\left\{2\int_{t_0}^s\frac{\alpha(u)}{T-u}\,\dd u\right\}\,\dd s\\
     & = C_1 + C_2 \int_{t_0}^t\exp\left\{2\int_{t_0}^s\frac{\alpha(u)}{T-u}\,\dd u\right\}\,\dd s, \nonumber
 \end{align}
 where
 \begin{align*}
   C_1 := \int_0^{t_0}\exp\left\{2\int_0^s\frac{\alpha(u)}{T-u}\,\dd u\right\}\,\dd s
   \qquad \text{and}\qquad
   C_2 := \exp\left\{2\int_0^{t_0}\frac{\alpha(u)}{T-u}\,\dd u\right\}.
 \end{align*}
Hence for all \ $t\in(t_0,T)$,
 \begin{align*}
   \langle M\rangle_t
     & \geq C_1 + C_2 \int_{t_0}^t\exp\left\{2\delta_1\int_{t_0}^s\frac{1}{T-u}\,\dd u\right\}\,\dd s
        =  C_1 + C_2 \int_{t_0}^t \left(\frac{T-t_0}{T-s}\right)^{2\delta_1} \,\dd s \\
     & = C_1 + C_2\frac{(T-t_0)^{2\delta_1}}{2\delta_1-1}
         \left( (T-t)^{1-2\delta_1} - (T-t_0)^{1-2\delta_1} \right),
 \end{align*}
 which yields that \ $\lim_{t\uparrow T} \langle M\rangle_t = \infty$, \
 since \ $\delta_1>1/2$.
\ Let us define the function \ $f:[1,\infty)\to (0,\infty)$ \ by \ $f(x):=x^{\delta_1/(2\delta_2-1)},$
 $x\geq 1.$
\ Then \ $f$ \ is strictly monotone increasing, since \ $\delta_1>0$ \ and \ $\delta_2>1/2$. \ Further
 \begin{equation}\label{fsquareint}
   \int_1^\infty \frac{1}{f(x)^2}\,\dd x
           = \int_1^\infty x^{-2\delta_1/(2\delta_2-1)}\,\dd x
           = \frac{2\delta_2-1}{1-2(\delta_2-\delta_1)}
           <\infty,
 \end{equation}
where we used that \ $2(\delta_2-\delta_1)-1<0$.
\ Hence we may apply a strong law of large numbers for continuous local martingales
 (see, e.g., $3^\circ)$ in Exercise 1.16 in Chapter V in Revuz and Yor \cite{RevYor}
 or Theorem 2.3 in Barczy and Pap \cite{BarPap2}) and then we obtain
 \begin{equation}\label{asconvx}
   \PP\left(\lim_{t\uparrow T} \frac{M_t}{f(\langle M\rangle_t)} = 0 \right)=1.
 \end{equation}
Further, for every \ $t\in[0,T)$ \ such that \ $\langle M\rangle_t\geq 1$ we have
\[
     X_t = \exp\left\{-\int_0^t\frac{\alpha(u)}{T-u}\,\dd u\right\} M_t
         = \exp\left\{-\int_0^t\frac{\alpha(u)}{T-u}\,\dd u\right\}
          f(\langle M\rangle_t)\,\frac{M_t}{f(\langle M\rangle_t)}.
 \]
Similarly as above, using \eqref{alphbounds}, \eqref{martres} and the monotonicity of \ $f$ \ we get
 \begin{align}\label{detpart}
  \exp&\left\{-\int_0^t\frac{\alpha(u)}{T-u}\,\dd u\right\} f(\langle M\rangle_t)\nonumber\\
      &\leq C_3\exp\left\{-\delta_1 \int_{t_0}^t\frac{1}{T-u}\,\dd u\right\}
       f\left(C_1 + C_2 \int_{t_0}^t\exp\left\{2\delta_2
                                \int_{t_0}^s\frac{1}{T-u}\,\dd u\right\}\,\dd s\right)\nonumber\\
     & = C_3\left(\frac{T-t}{T-t_0}\right)^{\delta_1}
         f\left(C_1 + C_2\frac{(T-t_0)^{2\delta_2}}{2\delta_2-1}
                        \big((T-t)^{1-2\delta_2} - (T-t_0)^{1-2\delta_2}\big)\right)\nonumber\\
     &\leq C_3\left(\frac{T-t}{T-t_0}\right)^{\delta_1}
          f\left(C_1 + C_2\frac{(T-t_0)^{2\delta_2}}{2\delta_2-1}
                                (T-t)^{1-2\delta_2}\right)\\
     & = C_3\left(\frac{T-t}{T-t_0}\right)^{\delta_1}
         \left(C_1 + C_2\frac{(T-t_0)^{2\delta_2}}{2\delta_2-1}
                                (T-t)^{1-2\delta_2}\right)^{\frac{\delta_1}{2\delta_2-1}}\nonumber\\
     & \to \frac{C_3}{(T-t_0)^{\delta_1}}
           \left(\frac{C_2(T-t_0)^{2\delta_2}}{2\delta_2-1}\right)^{\frac{\delta_1}{2\delta_2-1}}
           \quad \text{as} \quad t\uparrow T,\nonumber
 \end{align}
 where
 \[
   C_3:= \exp\left\{ - \int_0^{t_0} \frac{\alpha(u)}{T-u}\,\dd u\right\}.
 \]
Putting all together, we conclude that \ $\PP(\lim_{t\uparrow T}X_t=0)=1$ \ in case \ $\alpha(T)>1/2$.

Now we consider the case \ $0<\alpha(T)<1/2$. \
 Let \ $\delta_1$ \ and \ $\delta_2$ \ be given such that \ $0<\delta_1<\alpha(T)<\delta_2<\delta_1+1/2$
 \ and \ $\delta_2<1/2$.
\ Similarly as in the case \ $\alpha(T)>\frac{1}{2}$,
 \ from \eqref{alphbounds} and \eqref{martres} we get for all \ $t\in(t_0,T)$,
 \[
   \langle M\rangle_t
      \leq C_1 + C_2\frac{(T-t_0)^{2\delta_2}}{2\delta_2 - 1}
                    \left( (T-t)^{1-2\delta_2} - (T-t_0)^{1-2\delta_2} \right)
       \to C_1 + C_2\frac{(T-t_0)}{1-2\delta_2}
       \qquad \text{as \ $t\uparrow T$,}
 \]
 where we used that \ $\delta_2<1/2$.
\ This yields that \ $\lim_{t\uparrow T} \langle M\rangle_t<\infty$ \ if \ $\alpha(T)<1/2$
 \ (indeed, every bounded and monotone sequence is convergent).
Note also that for deriving \ $\lim_{t\uparrow T} \langle M\rangle_t<\infty$ \ we did not use
 that \ $\alpha(T)>0$, \ only that \ $\alpha(T)<1/2$.
\ Using Proposition 1.26 in Chapter IV and Proposition 1.8 in Chapter V in Revuz and Yor \cite{RevYor},
 we have the limit \ $M_T:=\lim_{t\uparrow T}M_t$ \ exists almost surely.
Since
 \[
   X_t = \exp\left\{-\int_0^t\frac{\alpha(u)}{T-u}\,\dd u\right\}M_t,\qquad t\in[0,T),
 \]
 and
 \[
  \exp\left\{-\int_0^{t}\frac{\alpha(u)}{T-u}\,\dd u\right\}
   \leq C_3\left(\frac{T-{t}}{T-t_0}\right)^{\delta_1}\to0
    \qquad \text{ as \ $t\uparrow T$,}
 \]
 we get \ $\PP(\lim_{t\uparrow T} X_t=0)=1$ \ also in case \ $0<\alpha(T)<1/2$.

 Finally, we consider the case \ $\alpha(T)=1/2$. \ Since the function
 \ $[0,T)\ni t\mapsto\langle M\rangle_t$ \
 is strictly increasing, \ we only have to consider the cases
 \ $\langle M\rangle_t\to\infty$ \
 or \ $\langle M\rangle_t\to\langle M\rangle_T :=\lim_{t\uparrow T} \langle M\rangle_t<\infty$
 \ as \ $t\uparrow T$.
\ If \ $\langle M\rangle_t\to\infty$, \ then \eqref{asconvx} and \eqref{detpart} are still valid,
 since \ $\delta_1>0$ \ and \ $1/2<\delta_2<\delta_1+1/2$, \ and hence $f$ is strictly increasing
 and \eqref{fsquareint} holds.
As in the case \ $\alpha(T)>1/2$ \ we conclude that \ $\PP(\lim_{t\uparrow T} X_t=0)=1$.
If \ $\langle M\rangle_t\to\langle M\rangle_T<\infty$, then \ $X_{t}\to 0$ \ almost surely
 as \ $t\uparrow T$ as in the case \ $\alpha(T)<1/2$.
\proofend

\begin{Rem}
If \ $\alpha(T)=\lim_{t\uparrow T}\alpha(t)$ \ exists and \ $\alpha(T)<0$, \ then
 there does not exist an almost surely continuous
 extension of the process \ $(X_t)_{t\in[0,T)}$ \ onto \ $[0,T]$ \ which would take some
 constant at time \ $T$ \ with probability one (i.e., which would be a bridge).
Indeed, the second moment of \ $X_t$ \ converges to infinity as \ $t\uparrow T$, \
 which can be checked as follows.
We get
 \begin{align*}
  \EE(X_t^2) = \exp\left\{-2\int_0^t\frac{\alpha(u)}{T-u}\,\dd u\right\}
               \EE(M_t^2)
              = \exp\left\{-2\int_0^t\frac{\alpha(u)}{T-u}\,\dd u\right\}
                \langle M\rangle_t,
                \qquad t\in[0,T),
 \end{align*}
 where \ $\lim_{t\uparrow T}\langle M\rangle_t<\infty$ \
 (as the proof of Theorem \ref{THEOREM8} shows) and for all \ $t\in[t_2,T)$,
 \ where \ $t_2$ \ is such that \ $\alpha(t)\leq \alpha(T)/2$, \ $t\in[t_2,T]$, \ we have
 \begin{align*}
 \exp\left\{-2\int_0^t\frac{\alpha(u)}{T-u}\,\dd u\right\}
   & \geq \exp\left\{-2\int_0^{t_2}\frac{\alpha(u)}{T-u}\,\dd u\right\}
         \exp\left\{-\alpha(T)\int_{t_2}^t\frac{1}{T-u}\,\dd u\right\} \\
   & = \exp\left\{-2\int_0^{t_2}\frac{\alpha(u)}{T-u}\,\dd u\right\}
       \left(\frac{T-t}{T-t_0}\right)^{\alpha(T)}
       \to\infty \qquad \text{as \ $t\uparrow T$.}
 \end{align*}
In case \ $\alpha(T)=0$ \ many things can happen concerning the limit behaviour of
 \ $X_t$ \ as \ $t\uparrow T$.
\ If  \ $\alpha$ \ is the identically zero function,
 then it is already argued in the Introduction that \ $X_T:=\lim_{t\uparrow T}X_t$ \ exists almost surely
 and has a nondegenerate Gauss distribution.
If \  $\alpha: [0,T)\to\RR$, $\alpha(t):=-(\log(T-t))^{-1}$, $t\in[0,T)$, \ and \ $t_0>T-1$, \ then
 for all \ $t\in[t_0,T)$,
 $$
 \int_{t_0}^t\frac{\alpha(u)}{T-u}\,\dd u=\log(-\log(T-t))-\log(-\log(T-t_0))\to\infty
 \qquad \text{as \ $t\uparrow T$,}
 $$
 which yields
 $$
 \exp\left\{ -\int_{0}^t\frac{\alpha(u)}{T-u}\,\dd u\right\}\to0.
 $$
As argued in the case \ $0< \alpha(T)<1/2$ \ of  the proof of Theorem \ref{THEOREM8},
 we have \ $M_T:=\lim_{t\uparrow T}M_t$ \ exists almost surely, hence \ $X_t\to0$ \ almost surely.
The same argument shows that if \ $\alpha:[0,T)\to\RR$, $\alpha(t):=(\log(T-t))^{-1}$,
  $t\in[0,T)$, \ then
 $$
  \exp\left\{ - \int_{0}^t\frac{\alpha(u)}{T-u}\,\dd u\right\}\to\infty.
  $$
 Using that \ $M_T$ \ is a non-degenerate normally distributed random variable with mean zero
 and variance \ $\lim_{t\uparrow T}\langle M\rangle_t$ \ (indeed, normally distributed random variables
 can converge in distribution only to a normally distributed random variable) we have \ $\PP(M_T=0)=0$,
 $\PP(M_T>0)=\PP(M_T<0)=1/2$ \ and hence
 \[
   \PP(\lim_{t\uparrow T} X_t =\infty) =  \PP(\lim_{t\uparrow T} X_t = -\infty) = \frac{1}{2}.
 \]
Especially, \ $X_t$ \ does not have a finite limit as \ $t\uparrow T$ \ almost surely.

 Finally, we remark that we do not know whether there exists an almost surely continuous extension
 in case the limit of \ $\alpha(t)$ \ as \ $t\uparrow T$ \ does not exist.
\proofend
\end{Rem}

Now we turn to the question of identical bridges for general $\alpha$-Wiener bridges.

 First we prove an auxiliary lemma (and a corollary of it) on the uniqueness of
 the drift and diffusion coefficients of the SDE \eqref{gen_OU_egyenlet}.
This result may be known but the authors were not able to find any reference for it.
We will only need part (iii) of the following lemma but the other parts may also
 be of independent interest.

\begin{Lem}\label{LEMMA3}
 Let \ $T>0$ \ be fixed and let us suppose that the processes \ $(Y^{(i)}_t)_{t\in[0,T)}$, $i=1,2$,
 are strong solutions of the SDEs
 \begin{align*}
   \begin{cases}
    \dd Y^{(i)}_t=b_i(t)\,Y^{(i)}_{t}\,\dd t + \sigma_i(t)\dd B_t^{(i)},\quad t\in[0,T),\\
    \phantom{\dd} Y^{(i)}_0 = \xi^{(i)},
   \end{cases}
   \qquad i=1,2,
 \end{align*}
 respectively, where \ $b_i$, $\sigma_i:[0,T)\to\RR$, $i=1,2,$ \ are
 continuous functions such that $\sigma_i(t)\ne 0$, $t\in[0,T)$, $i=1,2$,
 \ $(B_t^{(i)})_{t\geq 0},$ $i=1,2$, are one-dimensional standard Wiener processes
 and \ $\xi^{(i)}$, $i=1,2$, \ are Gauss random variables independent
 of \ $B^{(i)}$, $i=1,2$, \ respectively.
 \renewcommand{\labelenumi}{{\rm(\roman{enumi})}}
 \begin{enumerate}
   \item If the one-dimensional distributions of \ $Y^{(1)}$ \ and \ $Y^{(2)}$ \ coincide and
         $\EE\xi^{(1)} = \EE\xi^{(2)} \ne 0$, then \ $b_1(t)=b_2(t),$ $t\in[0,T)$, and
         \ $\vert \sigma_1(t)\vert = \vert \sigma_2(t)\vert$, $t\in[0,T)$.
   \item If the one-dimensional distributions of \ $Y^{(1)}$ \ and \ $Y^{(2)}$ \ coincide and
         $\sigma_1(t)=\sigma_2(t)=\sigma$, $t\in[0,T)$, \ for some \ $\sigma\in\RR$, $\sigma\ne0$,
         \ then \ $b_1(t)=b_2(t),$ $t\in[0,T)$.
   \item If the one- and two-dimensional distributions of \ $Y^{(1)}$ \ and \ $Y^{(2)}$ \ coincide,
         respectively, then \ $b_1(t)=b_2(t),$ $t\in[0,T)$, and
         \ $\vert \sigma_1(t)\vert = \vert \sigma_2(t)\vert$, $t\in[0,T)$.
 \end{enumerate}
\end{Lem}

\noindent{\bf Proof.}
 By Section 5.6 in Karatzas and Shreve \cite{KarShr}, we know that
 both SDEs have a unique strong solution that obeys an integral
 representation
 $$
    Y^{(i)}_t = \ee^{\overline b_i(t)}\xi^{(i)} +  \int_0^t\ee^{\overline b_i(t)-\overline b_i(s)}
                      \sigma_i(s)\,\dd B_s^{(i)},
                      \qquad t\in[0,T),\;\;\; i=1,2,
 $$
 where \ $\overline b_i(t):=\int_0^tb_i(s)\,\dd s,$ $t\in[0,T),$
 $i=1,2$, \ and strong uniqueness holds for both SDEs, see, e.g.,
 {\O}ksendal \cite[Theorem 5.2.1]{Oks} or Jacod and Shiryaev \cite[Chapter III, Theorems 2.32 and 2.33]{JacShi}.
By the assumptions, the one-dimensional marginals of the Gauss processes
 \ $ \ee^{\overline b_i(t)}\xi^{(i)} +  \int_0^t\ee^{\overline b_i(t)-\overline b_i(s)}
  \sigma_i(s)\,\dd B_s^{(i)},$ $t\in[0,T),$ $i=1,2,$ \ coincide.
Under the assumption that the one-dimensional distributions of \ $Y^{(1)}$ \ and
 \ $Y^{(2)}$ \ coincide, we have the means and the variances of these one-dimensional
 marginals are identical, namely,
 \begin{align}\label{SEGED4_suppl}
  \ee^{\overline b_1(t)} \EE\xi^{(1)}
      & = \ee^{\overline b_2(t)} \EE\xi^{(2)},\\ \label{SEGED4}
   \ee^{2\overline b_1(t)}\var(\xi^{(1)}) +
    \int_0^t\ee^{2(\overline b_1(t)-\overline b_1(s))}
       \sigma_1^2(s) \,\dd s
      & = \ee^{2\overline b_2(t)}\var(\xi^{(2)}) +
        \int_0^t\ee^{2(\overline b_2(t)-\overline b_2(s))}
         \sigma_2^2(s)  \,\dd s
 \end{align}
 for all \ $t\in[0,T)$.
Further, under the assumption that the one- and two-dimensional distributions
 of \ $Y^{(1)}$ \ and \ $Y^{(2)}$ \ coincide, respectively, besides \eqref{SEGED4_suppl}
 and \eqref{SEGED4} we also have the covariances of the coordinates of the two-dimensional
 marginals are identical, namely,
 \begin{align}\label{SEGED4_suppl2}
  \begin{split}
   \ee^{\overline b_1(s) + \overline b_1(t)}
     & \left[\var(\xi^{(1)}) + \int_0^{s\wedge t} \ee^{-2\overline b_1(u)}\sigma_1^2(u)\,\dd u\right]\\
     & =  \ee^{\overline b_2(s) + \overline b_2(t)}
       \left[\var(\xi^{(2)}) + \int_0^{s\wedge t} \ee^{-2\overline b_2(u)}\sigma_2^2(u)\,\dd u\right],
     \qquad s,t\in[0,T),
  \end{split}
 \end{align}
 see, e.g., Karatzas and Shreve \cite[(5.6.11)]{KarShr}.

(i): Let us suppose that the one-dimensional distributions
  of \ $Y^{(1)}$ \ and \ $Y^{(2)}$ \ coincide and \ $\EE\xi^{(1)} = \EE\xi^{(2)} \ne 0$.
\ By \eqref{SEGED4_suppl} we have
  \ $\ee^{\overline b_1(t)} = \ee^{\overline b_2(t)}$, $t\in[0,T)$, \
 and hence \ $\overline b_1(t) = \overline b_2(t)$, $t\in[0,T)$.
\ By differentiation with respect to \ $t$ \ and using also that \ $b_i$, $i=1,2$, \
 are continuous, we get \ $b_1(t) = b_2(t)$, $t\in[0,T)$.
Differentiating \eqref{SEGED4} with respect to \ $t$, we have
\begin{align*}
   & \ee^{2\overline b_1(t)} 2b_1(t)
        \left( \var(\xi^{(1)})
        + \int_0^t \ee^{-2\overline b_1(s)}\sigma_1^2(s) \,\dd s \right)
        +  \sigma_1^2(t)  \\
  & \phantom{\ee^{2\overline b_1(t)} 2b_1(t) \Big( \EE(\xi^{(1)})^2 +}
      = \ee^{2\overline b_2(t)} 2 b_2(t)
           \left( \var(\xi^{(2)})
         + \int_0^t \ee^{-2\overline b_2(s)} \sigma_2^2(s) \,\dd s \right)
         + \sigma_2^2(t), \quad t\in[0,T).
\end{align*}
By \eqref{SEGED4} and using also that we proved that the continuous functions
 \ $b_1$ \ and \ $b_2$ \ are equal, we get \ $\sigma_1^2(t)=\sigma_2^2(t),$ $t\in[0,T)$.

(ii): Let us suppose that the one-dimensional distributions
 of \ $Y^{(1)}$ \ and \ $Y^{(2)}$ \ coincide and
 \ $\sigma_1(t)=\sigma_2(t)=\sigma$, $t\in[0,T)$, \ for some \ $\sigma\in\RR$, $\sigma\ne 0$.
\ If \ $\EE\xi^{(1)} = \EE\xi^{(2)} \ne 0$, \ then the assertion follows by part (i).
If \ $\EE\xi^{(1)} = \EE\xi^{(2)} = 0$, \ then
differentiating \eqref{SEGED4} with respect to $t$, we have
 \begin{align*}
   & \ee^{2\overline b_1(t)} 2b_1(t)
        \left( \var(\xi^{(1)})
        + \sigma^2 \int_0^t \ee^{-2\overline b_1(s)}\,\dd s \right)
        +  \sigma^2  \\
  & \phantom{\ee^{2\overline b_1(t)} 2b_1(t) \Big( \EE(\xi^{(1)})^2 +}
      = \ee^{2\overline b_2(t)} 2 b_2(t)
           \left( \var(\xi^{(2)})
         + \sigma^2 \int_0^t \ee^{-2\overline b_2(s)} \,\dd s \right)
         + \sigma^2, \qquad t\in[0,T).
 \end{align*}
Using \eqref{SEGED4}, this yields that $b_1(t) = b_2(t)$, $t\in[0,T)$.

(iii): Let us suppose that the one- and two-dimensional distributions
 of \ $Y^{(1)}$ \ and \ $Y^{(2)}$ \ coincide, respectively.
For all fixed \ $s\in[0,T)$, \ differentiating \eqref{SEGED4_suppl2} with respect to
 \ $t$ \ on the interval \ $(s,T)$, \ we have
 \begin{align*}
    \ee^{\overline b_1(s) + \overline b_1(t)} b_1(t)
       & \left( \var(\xi^{(1)})
        + \int_0^s \ee^{-2\overline b_1(u)}\sigma_1^2(u) \,\dd u \right)  \\
   & = \ee^{\overline b_2(s) + \overline b_2(t)} b_2(t)
        \left( \var(\xi^{(2)})
        + \int_0^s \ee^{-2\overline b_2(u)}\sigma_2^2(u) \,\dd u \right),
        \qquad 0\leq s<t<T.
 \end{align*}
Then \eqref{SEGED4_suppl2} yields that \ $b_1(t) = b_2(t)$, $t\in(0,T)$, \ and
 the continuity of $b_1$ and $b_2$ implies that \ $b_1(0) = b_2(0)$. \
For all fixed \ $t\in(0,T)$, \ differentiating \eqref{SEGED4_suppl2} with respect to
 \ $s$ \ on the interval \ $(0,t)$, \ we have
 \begin{align*}
    \ee^{\overline b_1(s) + \overline b_1(t)} b_1(s)
       & \left( \var(\xi^{(1)})
        + \int_0^s \ee^{-2\overline b_1(u)}\sigma_1^2(u) \,\dd u \right)
        + \ee^{-\overline b_1(s) + \overline b_1(t)} \sigma_1^2(s)   \\
   & = \ee^{\overline b_2(s) + \overline b_2(t)} b_2(s)
        \left( \var(\xi^{(2)})
        + \int_0^s \ee^{-2\overline b_2(u)}\sigma_2^2(u) \,\dd u \right)
        + \ee^{-\overline b_2(s) + \overline b_2(t)} \sigma_2^2(s)
 \end{align*}
 for all \ $0<s<t<T$.
\ Since \ $b_1(t) = b_2(t)$, $t\in[0,T)$, \ was already checked,
 \eqref{SEGED4_suppl2} yields that \ $\sigma_1^2(t) = \sigma_2^2(t)$, $t\in(0,T)$, \ and the continuity of $\sigma_1$ and $\sigma_2$ implies that \ $\sigma_1^2(0) = \sigma_2^2(0)$.
\proofend

Next we formulate a simple corollary of Lemma \ref{LEMMA3}, which will be used several times later on
 in the proofs.

\begin{Cor}\label{CORROLARY1}
Let \ $T>0$ \ be fixed and let us suppose that the processes
 \ $(Y^{(i)}_t)_{t\in[0,T)}$, $i=1,2$, \ are
 strong solutions of the SDEs given in Lemma \ref{LEMMA3}.
Further, let \ $(\widetilde Y^{(i)}_t)_{t\in[0,T)}$, $i=1,2$, \ be almost surely
 continuous (Gauss) processes having the same finite dimensional distributions as
 \ $(Y^{(i)}_t)_{t\in[0,T)}$, $i=1,2$, \ respectively.
 \renewcommand{\labelenumi}{{\rm(\roman{enumi})}}
 \begin{enumerate}
   \item If the one-dimensional distributions of \ $\widetilde Y^{(1)}$ \ and
         \ $\widetilde Y^{(2)}$ \ coincide and $\EE\xi^{(1)} = \EE\xi^{(2)} \ne 0$, then
         \ $b_1(t)=b_2(t),$ $t\in[0,T)$, and \ $\vert \sigma_1(t)\vert = \vert \sigma_2(t)\vert$,
          $t\in[0,T)$.
   \item If the one-dimensional distributions of \ $\widetilde Y^{(1)}$ \ and \ $\widetilde Y^{(2)}$
         \ coincide and $\sigma_1(t)=\sigma_2(t)=\sigma$, $t\in[0,T)$, \  for some \ $\sigma\in\RR$,
         $\sigma\ne0$, \ then \ $b_1(t)=b_2(t),$ $t\in[0,T)$.
   \item If the one- and two-dimensional distributions of \ $\widetilde Y^{(1)}$ \ and
         \ $\widetilde Y^{(2)}$ \ coincide, respectively, then \ $b_1(t)=b_2(t),$ $t\in[0,T)$,
         and \ $\vert \sigma_1(t)\vert = \vert \sigma_2(t)\vert$, $t\in[0,T)$.
 \end{enumerate}
\end{Cor}

\noindent{\bf Proof.}
By the assumptions, the one-dimensional distributions of
 \ $Y^{(1)}$ \ and \ $Y^{(2)}$ \ coincide and hence Lemma \ref{LEMMA3} yields the assertion.
\proofend

\begin{Thm}\label{THEOREM4}
Let \ $T>0$ \ be fixed and \ $\alpha:[0,T)\to\RR$ \ be a continuous
function such that \ $\lim_{t\uparrow T}\alpha(t)\ne 1.$ \ There
does not exist an Ornstein-Uhlenbeck type process \ $(Z_t)_{t\geq 0}$
 \ given by the SDE \eqref{gen_OU_egyenlet} such that the law of the
 Ornstein-Uhlenbeck type bridge from \ $0$ \ to \ $0$ \ over the time-interval
 \ $[0,T]$ \ derived from \ $Z$ \ coincides with the law of the general $\alpha$-Wiener bridge.
\end{Thm}

\noindent{\bf Proof.}
We give an indirect proof. Let \ $(Z_t)_{t\geq 0}$ \ be an Ornstein-Uhlenbeck type
 process given by the SDE \eqref{gen_OU_egyenlet}.
Suppose that the law of the Ornstein-Uhlenbeck type bridge from \ $0$ \ to \ $0$ \ over
 \ $[0,T]$ \ derived from \ $Z$ \ coincides with the law of the general $\alpha$-Wiener bridge.
The process \ $(U_t)_{t\in[0,T)}$ \ given by \eqref{gen_OU_integral} with
 \ $a=0$ \ and \ $b=0$ \ is an Ornstein-Uhlenbeck type bridge from \ $0$ \ to \ $0$ \ over
 \ $[0,T]$ \ derived from \ $Z$.
\ By Lemma \ref{LEMMA4} and part (iii) of Corollary \ref{CORROLARY1}, we have
  \begin{align*}
    -\frac{\alpha(t)}{T-t}=q(t) - \sigma^2(t) \frac{\ee^{2(\bar q(T)-\bar q(t))}}{\gamma(t,T)},
        \quad t\in[0,T),\quad  \text{and}\quad|\sigma(t)|=1,\quad t\in[0,T).
  \end{align*}
Hence
 $$
    \alpha(t)=-(T-t)q(t)+\frac{(T-t)\ee^{2(\bar q(T)-\bar q(t))}}{\gamma(t,T)},
      \qquad t\in[0,T).
 $$
Using that \ $q$ \ is continuous, we have \ $\lim_{t\uparrow T}(T-t)q(t)=0\cdot q(T)=0,$ \
 and then
 \begin{align*}
  \lim_{t\uparrow T}\alpha(t)
     & =\lim_{t\uparrow T}\frac{(T-t)\ee^{2(\bar q(T)-\bar q(t))}}{\gamma(t,T)}
       =\lim_{t\uparrow T}\frac{-\ee^{2(\bar q(T)-\bar q(t))}-2q(t)(T-t)
        \ee^{2(\bar q(T)-\bar q(t))}}{-\sigma^2(t)\ee^{2(\bar q(T)-\bar q(t))}}\\
     &\, =1+2\lim_{t\uparrow T}q(t)(T-t)=1,
 \end{align*}
 where we used that
 \begin{equation}\label{gamma_derivative}
  \partial_1\gamma(u,T)=-\sigma^2(u)\ee^{2(\bar q(T)-\bar q(u))},\quad 0\leq u<T.
 \end{equation}
Hence we arrived at a contradiction.
 \proofend

The next remark shows that there exist general $\alpha$-Wiener bridges which are bridges
 derived from Ornstein-Uhlenbeck type processes.

\begin{Rem}\label{REMARK2}
Note that if \ $\alpha(t)=q(T-t)\coth(q(T-t)),$ $t\in[0,T),$ \ with some \ $q\not=0,$ \ then
  the SDE \eqref{gen_alpha_W_bridge} has the form
 \begin{align*}
  \begin{cases}
   \dd X_t=-q\coth(q(T-t))\,X_{t}\,\dd t+\dd B_t,\qquad t\in[0,T),\\
   \phantom{\dd} X_0=0,
  \end{cases}
 \end{align*}
 and, by Remark \ref{REMARK3}, this SDE coincides with the SDE satisfied by the
 Ornstein-Uhlenbeck bridge (given in \eqref{OU_hid1_intrep}) from $0$ to $0$ over $[0,T]$
 derived from the Ornstein-Uhlenbeck process given by the SDE
   \ $\dd Z_t=q\,Z_{t}\,\dd t+\dd B_t$, $t\geq 0$, \
 with an initial condition \ $Z_0$ \ having a Gauss distribution independent of the Wiener process \ $B$.
By L'Hospital's rule we have
\begin{align*}
   \lim_{t\uparrow T}\alpha(t)
   &  =\lim_{t\uparrow T}q(T-t)\coth(q(T-t))= \lim_{t\uparrow T}\frac{q(T-t)\cosh(q(T-t))}{\sinh(q(T-t))}\\
   &  =\lim_{t\uparrow T}\frac{-q\cosh(q(T-t))-q(T-t)\sinh(q(T-t))}{-q\cosh(q(T-t))}=1,
\end{align*}
as expected by Theorem \ref{THEOREM4}. This example shows that there are general
$\alpha$-Wiener bridges which can be derived from an
Ornstein-Uhlenbeck type process by taking a bridge.
For a more detailed discussion of this example, see Example \ref{REM_peldak1}.
\proofend
\end{Rem}

In what follows we will study the question whether
 every general $\alpha$-Wiener bridge with a continuously differentiable
 \ $\alpha$ \ such that \ $\lim_{t\uparrow T}\alpha(t)=1$ \ can be derived from
 some appropriate Ornstein-Uhlenbeck type process by taking a bridge, see Theorem \ref{THEOREM5}.
First, for our later purposes, we recall a well-known result about the solutions of special type
 of Riccati DEs, see, e.g., Reid \cite[Chapter I, Theorem 2.2]{Rei},
 Vrabie \cite[Theorems 1.3.4 and 1.3.5]{Vra} or Walter \cite[page 33]{Wal}.

\begin{Rem}\label{REMARK6}
Let \ $I\subset\RR$ \ be an interval, \ $I_0\subset I$ \ be a subinterval of \ $I$, \ $s\in I_0$, \ and
 \ $c:I\to\RR$ \ be a continuous function.
Further, let \ $w_0:I_0\to\RR$ \ be a solution of the Riccati type differential equation
 \begin{align}\label{Riccati_DE}
    w'(t) = -w^2(t) + c(t), \qquad t\in I_0.
 \end{align}
Then \ $w:I_0\to\RR$ \ is a solution of the DE \eqref{Riccati_DE} if and only if
 there exists a constant \ $C\in\RR\cup\{\infty\}$ \ such that
 \ $C\varphi(t)+\psi(t)\ne 0$, $t\in I_0$, \ and
 \[
   w(t) = w_0(t) + \frac{1}{C\varphi(t) + \psi(t)},\qquad t\in I_0,
 \]
 where, for \ $C\in\RR$,
 \ $u:=C\varphi + \psi$ \ is the unique solution of the DE
 \begin{align} \label{SEGED9}
   u'(t) - 2w_0(t)u(t) = 1, \qquad t\in I_0,
 \end{align}
 with an initial condition \ $u(s)=C$.
\ For \ $C=\infty$ \ we mean that \ $w(t)=w_0(t)$, $t\in I_0$.
\proofend
\end{Rem}

\begin{Thm}\label{THEOREM5}
Let \ $T>0$ \ be fixed and \ $\alpha:[0,T)\to\RR$ \ be a continuously differentiable
 function with \ $\lim_{t\uparrow T}\alpha(t)=1$.
 \begin{enumerate}
      \item[(i)] Let us consider the Ornstein-Uhlenbeck type process \ $(Z_t)_{t\geq 0}$ \ given by the SDE
            \eqref{gen_OU_egyenlet} with continuous functions  \ $q:[0,\infty)\to\RR$ \ and \ $\sigma:[0,\infty)\to\RR$
            and suppose that \ $q$ is continuously differentiable on \ $[0,T)$.
            If the law of the Ornstein-Uhlenbeck type bridge from \ $0$ \ to \ $0$ \
            over the time-interval \ $[0,T]$ \ derived from \ $Z$ \ coincides with the law of
            the general $\alpha$-Wiener bridge, then
            \begin{align}\label{gen_alpha_W_bridge_OU}
               q(t)=-\frac{\alpha(t)}{T-t}
                    + \frac{1}{C\exp\left\{-2\int_0^t\frac{\alpha(s)}{T-s}\,\dd s\right\}
                    +\int_0^t\exp\left\{-2\int_s^t\frac{\alpha(u)}{T-u}\,\dd u\right\}\,\dd s}
            \end{align}
            for all \ $t\in[0,T)$ \ with some \ $C\in(0,\infty)$, and
            \ $\sigma(t)=1,$ $t\in[0,T)$, \ or \ $\sigma(t)=-1,$ $t\in[0,T)$.
      \item[(ii)] Let \ $C\in(0,\infty)$ \ and define \ $q_C:[0,T)\to\RR$ \ as in
            \eqref{gen_alpha_W_bridge_OU}.
            If \ $\lim_{t\uparrow T}q_C(t)\in\RR$ \ exists,
            then there exists a continuous function
            \ $q:[0,\infty)\to\RR$ \ such that \ $q(t)=q_C(t)$, $t\in[0,T)$, \ and
            for all such extensions \ $q$, \ the law of the
            Ornstein-Uhlenbeck type bridge from \ $0$ \ to \ $0$ \ over the time-interval \ $[0,T]$ \
            derived from \ $Z$ \ given by the SDE \eqref{gen_OU_egyenlet} with
            \ $\sigma(t)=1,$ $t\in[0,T)$, \ or \ $\sigma(t)=-1,$ $t\in[0,T)$
            \ coincides with the law of the general $\alpha$-Wiener bridge.
 \end{enumerate}
\end{Thm}

\begin{Rem}\label{REMARK8_lacks}
Let \ $T>0$, \ $C>0$ \ and \ $\alpha:[0,T)\to\RR$ \ be a continuously differentiable function
  with \ $\lim_{t\uparrow T}\alpha(t)=1$.
\ We call the attention that if we define
 \ $q_C:[0,T)\to\RR$ \ as in \eqref{gen_alpha_W_bridge_OU}, then it is not sure that
 \ $\lim_{t\uparrow T}q_C(t)$ \ exists (see Example \ref{REM_peldak2}),
 which yields that it is not sure that \ $q_C$ \ can be continuously extended onto \ $[0,\infty)$.
\ Hence in this case the general $\alpha$-Wiener bridge can not be derived from an
 Ornstein-Uhlenbeck type process by taking a bridge, since in our setup the function
 \ $q$ \ in the SDE \eqref{gen_OU_egyenlet} should be defined on \ $[0,\infty)$.
 We also remark that the derivative of the function \ $\alpha$ \ does not appear
 in the formulation of our results in Theorem \ref{THEOREM5}, however we suppose that
 \ $\alpha$ \ is continuously differentiable.
The reason for this is our proof of technique, however
 one may get rid of this assumption using some other approach.
Finally, we emphasize that we were not able to derive a general sufficient condition on
 the function \ $\alpha$ \ such that in part (ii) of Theorem \ref{THEOREM5} the condition on
 the existence of the limit \ $\lim_{t\uparrow T}q_C(t)$ \ is satisfied.
 A special situation is discussed in Example \ref{REM_peldak2} in the next section.
\proofend
\end{Rem}

\noindent{\bf Proof of Theorem \ref{THEOREM5}.}
(i): Comparing the SDE \eqref{gen_alpha_W_bridge} with the SDE
 \eqref{gen_OU_hid1_egyenlet} for \ $a=0$ and $b=0,$ \
 part (iii) of Corollary \ref{CORROLARY1} implies that
 \begin{align}\label{gen_alpha_W_bridge_eq}
    -\frac{\alpha(t)}{T-t}=q(t) - \sigma^2(t)\frac{\ee^{2(\bar q(T)-\bar q(t))}}{\gamma(t,T)},
        \quad t\in[0,T),\quad \text{and}\quad|\sigma(t)|=1,\quad t\in[0,T).
  \end{align}
  It yields that \ $q(t)+\frac{\alpha(t)}{T-t}>0,$ $t\in[0,T),$ \ and hence
  $$
     \gamma(t,T)
       =\frac{\ee^{2(\bar q(T)-\bar q(t))}}{q(t)+\frac{\alpha(t)}{T-t}},
     \qquad t\in[0,T).
  $$
By differentiation with respect to \ $t$ \ and using \eqref{gamma_derivative},
we have for all \ $t\in[0,T)$,
 \begin{align*}
   -\ee^{2(\bar q(T)-\bar q(t))}
      =\frac{\ee^{2(\bar q(T)-\bar q(t))}(-2)q(t)\left(q(t)+\frac{\alpha(t)}{T-t}\right)
             -\ee^{2(\bar q(T)-\bar q(t))}\left(q'(t)+\frac{\alpha'(t)(T-t)+\alpha(t)}{(T-t)^2}\right)}
         {\left(q(t)+\frac{\alpha(t)}{T-t}\right)^2}.
 \end{align*}
Hence
 $$
  \left(q(t)+\frac{\alpha(t)}{T-t}\right)^2
     =q'(t)+2q(t)\left(q(t)+\frac{\alpha(t)}{T-t}\right)+\frac{\alpha'(t)(T-t)+\alpha(t)}{(T-t)^2},
      \quad t\in[0,T),
 $$
 which yields that
 \begin{align}\label{gen_alpha_W_bridge_DE}
   q'(t)=-q^2(t)+\frac{\alpha(t)(\alpha(t)-1)-\alpha'(t)(T-t)}{(T-t)^2},
       \qquad t\in[0,T).
 \end{align}
 Note that the differential equation \eqref{gen_alpha_W_bridge_DE} is of Riccati type.

By Remark \ref{REMARK6}, we get if \ $\widetilde q$ \ is a particular solution of the DE
 \eqref{gen_alpha_W_bridge_DE}, then the general solution of this  DE is
 $$
  \widetilde q+\frac{1}{C\varphi+\psi},\quad C\in\RR\cup\{+\infty\},
 $$
 where \ $u:=C\,\varphi+\psi,$ $C\in\RR,$ \ is a general solution of the inhomogeneous
  linear DE
 \begin{align}\label{gen_alpha_W_bridge_DE_inhom}
    u'(t)-2\widetilde q(t)u(t)=1,\qquad t\in[0,T),
 \end{align}
such that \ $u(t)\not=0,$  $t\in[0,T)$.
\ Now we check that \ $\widetilde q(t)=-\frac{\alpha(t)}{T-t},$ $t\in[0,T),$ \ is a solution of the DE
 \eqref{gen_alpha_W_bridge_DE}.
Indeed,
 $$
   \widetilde q'(t)=-\frac{\alpha'(t)(T-t)+\alpha(t)}{(T-t)^2},\qquad t\in [0,T),
 $$
 and
 \begin{align*}
  -\widetilde q^2(t)+\frac{\alpha(t)(\alpha(t)-1)-\alpha'(t)(T-t)}{(T-t)^2}
        &= -\frac{\alpha^2(t)}{(T-t)^2}
           + \frac{\alpha(t)(\alpha(t)-1)-\alpha'(t)(T-t)}{(T-t)^2}\\
        &= -\frac{\alpha'(t)(T-t) + \alpha(t)}{(T-t)^2},\qquad t\in[0,T).
 \end{align*}

\noindent The general solutions of the homogeneous linear DE
 \ $u'(t)-2\widetilde q(t)u(t)=0$, $t\in[0,T)$, \ which corresponds to the inhomogeneous
 linear DE \eqref{gen_alpha_W_bridge_DE_inhom} are
 \begin{align*}
   u(t)=C\exp\left\{\int_0^t2 \widetilde q(s)\,\dd s\right\}
        =C\exp\left\{-2\int_0^t\frac{\alpha(s)}{T-s}\,\dd s\right\},
         \quad t\in[0,T),\quad C\in\RR.
 \end{align*}
 Now we are searching for a particular solution of the DE
 \eqref{gen_alpha_W_bridge_DE_inhom} by the method of variation of
 constants. Let
 $$
   u_0(t):=c(t)\exp\left\{-2\int_0^t\frac{\alpha(s)}{T-s}\,\dd s\right\},\quad t\in[0,T),
 $$
  be a (particular) solution of the DE \eqref{gen_alpha_W_bridge_DE_inhom}.
 Then
\begin{align*}
  1&=u'_0(t)-2\widetilde q(t)u_0(t)\\
    &=c'(t)\exp\left\{-2\int_0^t\frac{\alpha(s)}{T-s}\,\dd s\right\}
     +c(t)\exp\left\{-2\int_0^t\frac{\alpha(s)}{T-s}\,\dd s\right\}
       (-2)\frac{\alpha(t)}{T-t}\\
    &\phantom{=\;}
     +2\frac{\alpha(t)}{T-t}c(t)\exp\left\{-2\int_0^t\frac{\alpha(s)}{T-s}\,\dd s\right\},
     \quad t\in[0,T),
\end{align*}
which yields that
$$
  c'(t)=\exp\left\{2\int_0^t\frac{\alpha(s)}{T-s}\,\dd s\right\},
   \quad t\in[0,T),
$$
and hence we may choose
$$
   c(t)=\int_0^t\exp\left\{2\int_0^s\frac{\alpha(u)}{T-u}\,\dd u\right\}\,\dd s,
        \quad t\in[0,T),
$$
  and in this case
$$
   u_0(t)=\int_0^t\exp\left\{-2\int_s^t\frac{\alpha(u)}{T-u}\,\dd u\right\}\,\dd s,
           \quad t\in[0,T).
$$
Hence the general solution of the DE
\eqref{gen_alpha_W_bridge_DE_inhom} is
$$
  u(t)=C\exp\left\{-2\int_0^t\frac{\alpha(s)}{T-s}\,\dd s\right\}
       +\int_0^t\exp\left\{-2\int_s^t\frac{\alpha(u)}{T-u}\,\dd u\right\}\,\dd s,
       \quad t\in[0,T),\;\;C\in\RR.
$$
This yields that the general solution of the DE
\eqref{gen_alpha_W_bridge_DE} is
$$
  q_C(t)=-\frac{\alpha(t)}{T-t}
            +\frac{1}{C\exp\left\{-2\int_0^t\frac{\alpha(s)}{T-s}\,\dd s\right\}
            +\int_0^t\exp\left\{-2\int_s^t\frac{\alpha(u)}{T-u}\,\dd u\right\}\,\dd s},
   \quad t\in[0,T),
$$
where \ $C\in\RR\cup\{+\infty\}$ \ is such that the denominator
$$
 C\exp\left\{-2\int_0^t\frac{\alpha(s)}{T-s}\,\dd s\right\}
            +\int_0^t\exp\left\{-2\int_s^t\frac{\alpha(u)}{T-u}\,\dd u\right\}\,\dd s
$$
 is not \ $0$ \ for any \ $t\in[0,T).$ \
With \ $C=+\infty$ \ we get (back) the solution \
 $\widetilde q(t)=-\frac{\alpha(t)}{T-t},$ $t\in[0,T)$, and for \ $C\in\RR$ \
 we have
$$
  C\exp\left\{-2\int_0^t\frac{\alpha(s)}{T-s}\,\dd s\right\}
            +\int_0^t\exp\left\{-2\int_s^t\frac{\alpha(u)}{T-u}\,\dd u\right\}\,\dd s
      \not =0,\qquad t\in[0,T),
$$
if and only if
 $$
   \int_0^t\exp\left\{2\int_0^s\frac{\alpha(u)}{T-u}\,\dd u\right\}\,\dd s
      \not= -C,
      \quad t\in[0,T).
 $$
This implies that the general solution of the DE \eqref{gen_alpha_W_bridge_DE}
 is \ $q_C(t),$ $t\in[0,T),$ \ where \ $C\in\RR\cup\{+\infty\}$ \
 is such that
 \begin{align*}
   C\not\in\left(-\int_0^T\exp\left\{2\int_0^t\frac{\alpha(u)}{T-u}\,\dd u\right\}\,\dd t,0\right].
 \end{align*}
 By the proof of Theorem \ref{THEOREM8}, we have \ $\lim_{t\uparrow T}\alpha(t)=1$ \ yields that
 \[
 \int_0^T\exp\left\{2\int_0^t\frac{\alpha(u)}{T-u}\,\dd u\right\}\,\dd t =\infty,
 \]
 and hence \ $C\in(0,\infty]$.

Hence the general solution of the equation \eqref{gen_alpha_W_bridge_eq} is
 \ $q_C(t),$ \ $t\in[0,T),$ \ where \ $C\in(0,\infty)$. \
Indeed, the case \ $C=+\infty$ \ has to be excluded, since
 with \ $C=+\infty,$ \ $q_C(t)=-\frac{\alpha(t)}{T-t},$
 $t\in[0,T),$ \ and in this case it does not hold that
 \ $q_C(t)+\frac{\alpha(t)}{T-t}>0$, $t\in[0,T)$
 (which should be valid, see the beginning of the proof),
 further, by the assumption \ $\lim_{t\uparrow T}\alpha(t)=1$,
 \ we have \ $\lim_{t\uparrow T}q_C(t)=-\infty$, \ which yields that \ $q_C$ \ can not be
 extended to a continuous function onto \ $[0,\infty)$.
We give another brief explanation why we have to exclude the case \ $C=+\infty.$
\ With \ $C=+\infty$ \ we have \ $q_C(t)=-\frac{\alpha(t)}{T-t},$ $t\in[0,T),$
 \ and thus we would already start with the SDE \eqref{gen_alpha_W_bridge} of the general
 $\alpha$-Wiener bridge.
Hence we would try to derive a bridge from the bridge itself, which is not allowed with
 the procedure described in Section \ref{Section_OU_proc_bridge}.
 Indeed, for the transition densities of the bridge (see formula \eqref{gen_bridge_densities}) we need to know
 the transition density \ $p_{s,T}^{Z}(x,0)$ \ for the bridge itself which is not defined.

(ii): The possibility of such an extension follows readily.
Comparing the integral representation \eqref{gen_alpha_W_bridge_intrep}
 of the general $\alpha$-Wiener bridge with the integral
 representation \eqref{gen_OU_integral} of the Ornstein-Uhlenbeck type
 bridge for \ $a=0$ and $b=0$, \ by the definition of the general \ $\alpha$-Wiener bridge
 and Definition \ref{DEFINITION_bridge},
 it is enough to check that
 $$
   \exp\left\{-\int_s^t\frac{\alpha(u)}{T-u}\,\dd u\right\}
      =\frac{\gamma_C(t,T)}{\gamma_C(s,T)} \ee^{\bar q_C(t)-\bar q_C(s)},
      \qquad 0\leq s\leq t<T,
 $$
 where
 \[
    \gamma_C(s,t) = \int_s^t \ee^{2(\bar q_C(t)-\bar q_C(u))}\,\dd u,
                        \qquad 0\leq s<t.
 \]
 Indeed, for the case \ $\sigma(t)=-1$, $t\in[0,T)$, \ we note that if we replace
 the Wiener process \ $B$ \ with \ $-B$ \ in \eqref{gen_OU_integral} we still have an
 integral representation of the Ornstein-Uhlenbeck type bridge from \ $0$ \ to \ $0$
 \ over \ $[0,T]$ \ derived from \ $Z$.
\ Using that the function \ $q_C$ \ satisfies the equation \eqref{gen_alpha_W_bridge_eq},
 by \eqref{gamma_derivative}, we get for all \ $0\leq s\leq t<T,$ \
 \begin{align*}
  \exp\left\{-\int_s^t\frac{\alpha(u)}{T-u}\,\dd u\right\}
    & = \exp\left\{\int_s^t
         \left(q_C(u)-\frac{\ee^{2(\bar q_C(T)-\bar q_C(u))}}{\gamma(u,T)}\right)\,\dd u\right\}\\
    & = \exp\Big\{\bar q_C(t)-\bar q_C(s) + \ln(\gamma_C(t,T)) - \ln(\gamma_C(s,T)) \Big\}\\
    & = \ee^{\bar q_C(t)-\bar q_C(s)}\frac{\gamma_C(t,T)}{\gamma_C(s,T)},
 \end{align*}
as desired.
 \proofend

\begin{Rem}
 Note that in Theorem \ref{THEOREM5} the condition
 \ $\lim_{t\uparrow T}\alpha(t)=1$ \ on the function
 \ $\alpha$ \ is necessary. For this, note that a function \ $q$ \
 given in Theorem \ref{THEOREM5} satisfies the equation \eqref{gen_alpha_W_bridge_eq}
 and hence, by the proof of Theorem \ref{THEOREM4}, we have
 \ $\lim_{t\uparrow T}\alpha(t)=1.$
\proofend
\end{Rem}

\section{ Examples }\label{Section_spec_alpha_Wiener}

First we give examples of continuous functions \ $\alpha:[0,T)\to\RR$ \ such that
 \ $\lim_{t\uparrow T}\alpha(t)=1$ \ and, depending on a parameter, the general $\alpha$-Wiener bridge
 either can not be derived from any Ornstein-Uhlenbeck type process, or it can be derived from
 an Ornstein-Uhlenbeck type process by taking a bridge.

\begin{Ex}\label{REM_peldak2}
Let \ $T>0$ \ and \ $\alpha:[0,T)\to\RR$, $\alpha(t):=1\pm(T-t)^\beta$, $t\in[0,T)$,
 \ for some \ $\beta>0$.
\ For \ $0\leq s<t<T$ \ we have
 \begin{align*}
  \int_s^t\frac{\alpha(u)}{T-u}\,\dd u &
   = \int_s^t \frac{1\pm(T-u)^\beta}{T-u}\,\dd u
   = \int_s^t \left(\frac{1}{T-u}\pm(T-u)^{\beta-1} \right)\,\dd u\\
   & = \ln\left(\frac{T-s}{T-t}\right) \mp \frac1\beta\left((T-t)^\beta-(T-s)^\beta\right),
  \end{align*}
 and especially
 \[
  \int_0^t\frac{\alpha(u)}{T-u}\,\dd u
  = \ln\left(\frac{T}{T-t}\right)\mp \frac1\beta\left((T-t)^\beta-T^\beta\right)
    \to \infty \qquad \text{as \ $t\uparrow T$.}
 \]
Hence
 \begin{align*}
  \exp\left\{-2\int_0^t\frac{\alpha(u)}{T-u}\,\dd u \right\}
             = \left(\frac{T-t}{T}\right)^2\exp\left\{
                      \pm\frac2\beta\left((T-t)^\beta-T^\beta\right)\right\}
             \to 0\qquad \text{as \ $t\uparrow T$}
 \end{align*}
 and
 \begin{align*}
  \int_0^t\exp\left\{-2\int_s^t\frac{\alpha(u)}{T-u}\,\dd u \right\}\,\dd s
  & = \int_0^t\left(\frac{T-t}{T-s}\right)^2
    \exp\left\{ \pm\frac2\beta \left((T-t)^\beta-(T-s)^\beta\right)\right\}\,\dd s\\
  & = (T-t)^2\exp\left\{ \pm\frac2\beta(T-t)^\beta\right\}
      \int_0^t\frac{\exp\left\{\mp\frac2\beta(T-s)^\beta\right\}}{(T-s)^2}\,\dd s.
 \end{align*}
Then for all \ $C\in(0,\infty)$ \ the function \ $q_C$ \ in Theorem \ref{THEOREM5}
 takes the form
  \begin{align*}
   q_C(t) & = - \frac{1 \pm (T-t)^\beta}{T-t}
          +\Bigg(C\left(\frac{T-t}{T}\right)^2
                   \exp\left\{\pm\frac2\beta\left((T-t)^\beta-T^\beta\right)\right\}\\
          & \phantom{=} + (T-t)^2\exp\left\{\pm\frac2\beta\left((T-t)^\beta-T^\beta\right)\right\}
             \int_0^t \frac{\exp\left\{\mp\frac2\beta\left((T-s)^\beta-T^\beta\right)\right\}}
                           { (T-s)^2}\,\dd s \Bigg)^{-1} \\
          & =  -\frac{\frac{C}{T^2}\left(\frac1{T-t}\pm(T-t)^{\beta-1}\right)}
                {\frac{C}{T^2}
                 +\int_0^t \frac{\exp\left\{ \mp\frac2\beta\left((T-s)^\beta-T^\beta\right)\right\}}
                   {(T-s)^2}\,\dd s}\\
          & \phantom{=\;}
               - \frac{\left(1 \mp(T-t)^\beta\right)\int_0^t \frac{\exp\left\{\mp
                  \frac2\beta\left((T-s)^\beta-T^\beta\right)\right\}}{(T-s)^2}\,\dd s
                  -\frac{\exp\left\{\mp\frac2\beta\left((T-t)^\beta-T^\beta\right)\right\}}{T-t}}
                   {\frac{C(T-t)}{T^2}+(T-t)\int_0^t
                   \frac{\exp\left\{\mp
                    \frac2\beta\left((T-s)^\beta-T^\beta\right)\right\}}{(T-s)^2}\,\dd s}
 \end{align*}
Here, we get
 \begin{align*}
   &\int_0^t \frac{\exp\left\{\mp\frac2\beta\left((T-s)^\beta-T^\beta\right)\right\}}{(T-s)^2}\,\dd s
      \geq \left(\min_{u\in[0,T]} \exp\left\{\mp\frac2\beta
         \left((T-u)^\beta-T^\beta\right)\right\} \right)
         \int_0^t \frac{1}{(T-s)^2}\,\dd s\\
   &\phantom{\qquad\qquad}
      = \left(\min_{u\in[0,T]} \exp\left\{\mp
         \frac2\beta\left((T-u)^\beta-T^\beta\right)\right\}\right) \left( \frac{1}{T-t} - \frac{1}{T}\right)
    \to\infty \qquad \text{as \ $t\uparrow T,$}
 \end{align*}
 and
 $$
  \frac1{T-t}\pm(T-t)^{\beta-1} = \frac{1\pm (T-t)^\beta}{T-t}\to\infty
      \qquad \text{as \ $t\uparrow T,$}
 $$
 hence, by L'Hospital's rule we have
 \begin{align*}
  \lim_{t\uparrow T}
    \frac{\frac{C}{T^2}\left(\frac1{T-t} \pm(T-t)^{\beta-1}\right)}{\frac{C}{T^2}
      +\int_0^t \frac{\exp\left\{ \mp\frac2\beta\left((T-s)^\beta-T^\beta\right)\right\}}{(T-s)^2}\,\dd s}
  & = \lim_{t\uparrow T}
     \frac{\frac{C}{T^2}\left(\frac1{(T-t)^2}
       \mp(\beta-1)(T-t)^{\beta-2}\right)}
       {\frac{\exp\left\{ \mp\frac2\beta\left((T-t)^\beta-T^\beta\right)\right\}}{(T-t)^2}} \\
     &   = \frac{C}{T^2}\exp\left\{ \mp\frac2\beta\,T^\beta\right\}
 \end{align*}
 and
 \begin{align*}
  \lim_{t\uparrow T} (T-t)  \int_0^t \frac{\exp\left\{\mp
      \frac2\beta\left((T-s)^\beta-T^\beta\right)\right\}}{(T-s)^2}\,\dd s
    & = \lim_{t\uparrow T}\frac{\frac{\exp\left\{\mp
        \frac2\beta\left((T-t)^\beta-T^\beta\right)\right\}}{(T-t)^2}}{\frac1{(T-t)^2}}
      = \exp\left\{  \pm\frac2\beta \,T^\beta\right\}.
 \end{align*}
  Then we get
 \begin{align*}
  \lim_{t\uparrow T} q_C(t)
   & = - \lim_{t\uparrow T} \left\{\left(1 \mp(T-t)^\beta\right)
         \int_0^t \frac{\exp\left\{ \mp\frac2\beta (T-s)^\beta \right\}}{(T-s)^2}\,\dd s
      -\frac{\exp\left\{ \mp\frac2\beta  (T-t)^\beta \right\}}{T-t}\right\}\\
   & \phantom{=\;\,} -\frac{C}{T^2}\exp\left\{ \mp\frac2\beta\,T^\beta\right\}.
 \end{align*}
Since integration by parts yields that for \ $t\in(0,T)$
 \begin{align*}
 \int_0^t \frac{\exp\left\{\mp\frac2\beta (T-s)^\beta \right\}}{(T-s)^2}\,\dd s
  & = \frac{\exp\left\{ \mp\frac2\beta (T-t)^\beta \right\}}{T-t}
      -\frac{\exp\left\{ \mp\frac2\beta\,T^\beta\right\}}{T}\\
  & \phantom{=\;} \mp 2\int_0^t \frac{\exp\left\{ \mp
     \frac2\beta (T-s)^\beta \right\}}{(T-s)^{2-\beta}}\,\dd s,
 \end{align*}
 we have
 \begin{align*}
  \lim_{t\uparrow T} q_C(t)
   = \left(-\frac{C}{T^2} + \frac{1}{T}\right)\exp\left\{ \mp\frac2\beta T^\beta\right\}
   &  \pm2\lim_{t\uparrow T}\int_0^t
      \frac{\exp\left\{ \mp\frac2\beta(T-s)^\beta \right\}}
            {(T-s)^{2-\beta}}\,\dd s\\
   & \pm \lim_{t\uparrow T}(T-t)^\beta
      \int_0^t \frac{\exp\left\{\mp\frac2\beta(T-s)^\beta \right\}}{(T-s)^{2}}\,\dd s.
 \end{align*}
The function \ $[0,T]\ni t\mapsto\exp\left\{ \mp\frac2\beta(T-t)^\beta\right\}$ \ is bounded
 and hence the first limit above exists if and only if \ $\int_0^T\frac{1}{(T-s)^{2-\beta}}\,\dd s<\infty$,
 \ thus  if and only if \ $\beta>1$.
\ On the other hand, for the second limit above we get by L'Hospital's rule
 \begin{align*}
 \lim_{t\uparrow T}(T-t)^\beta\int_0^t
 \frac{\exp\left\{\mp\frac2\beta(T-s)^\beta\right\}}{(T-s)^{2}}\,\dd s
 & =  \lim_{t\uparrow T}
    \frac{\frac{\exp\left\{\mp\frac2\beta(T-t)^\beta\right\}}{(T-t)^2}}
                  {\beta(T-t)^{-\beta-1}}
  = \frac1\beta \lim_{t\uparrow T}(T-t)^{\beta-1},
 \end{align*}
 which exists if and only if \ $\beta\geq1$.

Alltogether we conclude that for \ $\alpha(t)=1\pm (T-t)^\beta$, $t\in[0,T)$,
 \ for some \ $\beta>0$, \ the limit of \ $q_C(t)$ \ as \ $t\uparrow T$ \
 exists if and only if \ $\beta>1$,
 i.e., for the given function \ $\alpha$, \ the general \ $\alpha$-Wiener bridge can be
 derived from an Ornstein-Uhlenbeck type bridge by taking a bridge
 if and only if \ $\beta>1$.

Further, we note that for the given function \ $\alpha$, \
 the limit of the ''inhomogeneity part'' of the Riccati type DE \eqref{gen_alpha_W_bridge_DE} is
 \begin{align*}
   \lim_{t\uparrow T} \frac{\alpha(t)(\alpha(t)-1)-\alpha'(t)(T-t)}{(T-t)^2}
   &   = \lim_{t\uparrow T} \frac{\pm(1 \pm(T-t)^\beta)(T-t)^\beta\pm
          \beta(T-t)^\beta}{(T-t)^2}\\
    &  = \pm\lim_{t\uparrow T}(T-t)^{\beta-2}\left(1+\beta \pm(T-t)^\beta\right)
 \end{align*}
 and this limit exists if and only if \ $\beta\geq2$.
\proofend
\end{Ex}

 Next, using Theorem \ref{THEOREM5}, we give a detailed study of the example presented
 in Remark \ref{REMARK2}.

\begin{Ex}\label{REM_peldak1}
Let $T>0$, $q\ne 0$, and $\alpha:[0,T)\to\RR$,
 $\alpha(t):=q(T-t)\coth(q(T-t))$, $t\in[0,T)$.
In Remark \ref{REMARK2}, without using Theorem \ref{THEOREM5}, we checked that the law of the general
 $\alpha$-Wiener bridge (with the given function \ $\alpha$) \ coincides with the law of the
 Ornstein-Uhlenbeck bridge from \ $0$ \ to \ $0$ \ over \ $[0,T]$ \ derived from the Ornstein-Uhlenbeck
 process given by the SDE
 \[
    \dd Z_t = q Z_t\,\dd t + \dd B_t,\quad t\geq 0,
 \]
 with an initial condition \ $Z_0$ \ having a Gauss distribution independent of the Wiener
 process \ $B$.
\ In what follows we give a presentation of this result using Theorem \ref{THEOREM5}
 in order to give an application of this theorem.
Namely, let
 \[
    C:=\frac{1}{q(1+\coth(qT))}>0,
 \]
 and let us define the function \ $q_C:[0,T)\to\RR$ \ as in \eqref{gen_alpha_W_bridge_OU}.
Since for all \ $0\leq s<t<T$,
 \[
    \int_s^t \frac{\alpha(u)}{T-u}\,\dd u = q \int_s^t \coth(q(T-u)) \,\dd u
                                          = -\frac{1}{q}\ln\left(\frac{\sinh(q(T-t))}{\sinh(q(T-s))}\right),
 \]
 we have for all \ $t\in[0,T)$,
 \begin{align*}
   q_C(t) & = -q\coth(q(T-t))
             + \frac{1}{C\left(\frac{\sinh(q(T-t))}{\sinh(qT)}\right)^2
                        + \DS\int_0^t \left(\frac{\sinh(q(T-t))}{\sinh(q(T-s))}\right)^2 \,\dd s}\\
          & = -q\coth(q(T-t))
              + \frac{1}{C\left(\frac{\sinh(q(T-t))}{\sinh(qT)}\right)^2
                        + \frac{(\sinh(q(T-t)))^2}{q}
                            \big(\coth(q(T-t)) - \coth(qT)\big)}\\
          & = -q\coth(q(T-t))
              + \frac{q(\sinh(q(T-t)))^{-2}}
                  {(1+\coth(qT))^{-1}(\sinh(qT))^{-2} + \coth(q(T-t)) - \coth(qT) }.
 \end{align*}
Since for all \ $x\ne 0$,
\begin{align*}
   (1+\coth(x))^{-1}&(\sinh(x))^{-2} - \coth(x)
      = \frac{1}{(\sinh(x))^2 + \cosh(x)\sinh(x)}
        - \frac{\cosh(x)}{\sinh(x)}\\
    & = \frac{1- \cosh(x)\sinh(x) - (\cosh(x))^2}{\sinh(x)(\sinh(x)+\cosh(x))}
     = \frac{-(\sinh(x))^2 - \cosh(x)\sinh(x)}{\sinh(x)(\sinh(x)+\cosh(x))}
      = -1,
 \end{align*}
 we have
 \[
   q_C(t) = -q\coth(q(T-t)) + \frac{q(\sinh(q(T-t)))^{-2}}{\coth(q(T-t)) -1 },
    \qquad t\in[0,T).
 \]
Since for all \ $x\ne 0$,
\begin{align*}
   -\coth(x) & + \frac{(\sinh(x))^{-2}}{\coth(x)-1}
      = \frac{1}{\sinh(x)}\left(-\cosh(x) + \frac{1}{\cosh(x) - \sinh(x)} \right)\\
    & = \frac{-(\cosh(x))^2 + \cosh(x)\sinh(x) +1 }{\sinh(x)(\cosh(x) - \sinh(x))}
      = \frac{-(\sinh(x))^2 + \cosh(x)\sinh(x)}{\sinh(x)(\cosh(x) - \sinh(x))}
      =1,
 \end{align*}
 we have \ $q_C(t)=q$, $t\in[0,T)$.
\ Hence \ $\lim_{t\uparrow T}q_C(t)=q$, \ and then part (ii) of Theorem \ref{THEOREM5}
 yields the desired statement (formulated at the beginning of the example).
\proofend
\end{Ex}

In what follows we examine the special case of general $\alpha$-Wiener bridges
 with constant $\alpha\in\RR$, i.e., we specialize our results for (usual) $\alpha$-Wiener bridges.
As an immediate consequence of Theorem \ref{THEOREM4} we get:

\begin{Cor}\label{THEOREM3}
Let \ $T>0$ \ and \ $\alpha\in\RR$ \ be fixed such that \ $\alpha\not= 1.$ \
There does not exist an Ornstein-Uhlenbeck type process \ $(Z_t)_{t\geq 0}$
 \ given by the SDE \eqref{gen_OU_egyenlet} such that the law of
 the Ornstein-Uhlenbeck type bridge from \ $0$ \ to \ $0$ \ over \ $[0,T]$ \ derived
 from \ $Z$ \ coincides with the law of the $\alpha$-Wiener bridge.
\end{Cor}

\begin{Rem}\label{REMARK7}
We note that in case of \ $\alpha\leq 0$, \ the assertion of Corollary \ref{THEOREM3}
 follows immediately without reference to Theorem \ref{THEOREM4}.
Indeed, Ornstein-Uhlenbeck type bridges on the time-interval \ $[0,T]$ \ are almost surely
 constant at time \ $T$, \ however in case of \ $\alpha\leq 0$ \ this property does not
 hold for an \ $\alpha$-Wiener bridge as it was detailed in the introduction.
\proofend
\end{Rem}

The next theorem is a special case of Theorem \ref{THEOREM5} for (usual) $\alpha$-Wiener bridges.

\begin{Thm}\label{THEOREM6}
 Let \ $T>0$ \ be fixed. Let us consider the Ornstein-Uhlenbeck type process
 \ $(Z_t)_{t\geq 0}$ \ given by the SDE \eqref{gen_OU_egyenlet}
 with continuous functions  \ $q:[0,\infty)\to\RR$ \ and \ $\sigma:[0,\infty)\to\RR$
 and suppose that $q$ is continuously differentiable on \ $[0,T)$.
\ Then the law of the Ornstein-Uhlenbeck type bridge from \ $0$ \ to \ $0$ \ over
 \ $[0,T]$ \ derived from \ $Z$ \ coincides with the law of the usual Wiener
 bridge from $0$ to $0$ over $[0,T]$ if and only if
$$
  q(t)=\frac{1}{t+C},\qquad t\in[0,T),
$$
 with some \ $C\in(\RR\setminus[-T,0])\cup\{+\infty\},$ \ and
 \  $\sigma(t)=1$, $t\in[0,T)$, \ or \ $\sigma(t)=-1$, $t\in[0,T)$.
\ Note that for \ $C=\infty$ \ we mean that \ $q(t)=0$, $t\in[0,T)$.
\end{Thm}

\noindent{\bf Proof.}
We check that Theorem \ref{THEOREM5} implies Theorem \ref{THEOREM6}.
By assumption the conditions of part (i) of Theorem \ref{THEOREM5} are
 satisfied with \ $\alpha(t) :=1,$ $t\in[0,T).$ \ Hence
 the set of continuous functions \ $q:[0,\infty)\to\RR$ \ which are continuously differentiable
 on \ $[0,T)$ \ and for which the law of the Ornstein-Uhlenbeck type bridge
 from $0$ to $0$ over \ $[0,T]$ \ derived from \ $Z$ \
 coincides with the law of the $\alpha$-Wiener bridge with \ $\alpha=1$ \
 (i.e., the usual Wiener bridge from $0$ to $0$ over \ $[0,T]$) \
 can be parametrized as \ $q_C,$ $C>0,$ \ where for all \ $C>0$ \
$$
  q_C(t)
     =-\frac{1}{T-t}
       +\frac{1}{C\exp\left\{-2\int_0^t\frac{1}{T-s}\,\dd s\right\}
                 +\int_0^t\exp\left\{-2\int_s^t\frac{1}{T-u}\,\dd u\right\}\,\dd s},
        \qquad t\in[0,T).
$$
Then
\begin{align*}
 q_C(t)
    &=-\frac{1}{T-t}
      +\frac{1}{C\exp\big\{2(\ln(T-t)-\ln(T))\big\}
         +\int_0^t\exp\big\{2(\ln(T-t)-\ln(T-s))\big\}\,\dd s}\\
    &=-\frac{1}{T-t}
      +\frac{1}{C\left(\frac{T-t}{T}\right)^2
               +\int_0^t\left(\frac{T-t}{T-s}\right)^2\,\dd s}
     =-\frac{1}{T-t}
        +\frac{1}{C\left(\frac{T-t}{T}\right)^2
          +(T-t)^2\left(\frac{1}{T-t}-\frac{1}{T}\right)}\\
    &=-\frac{1}{T-t}
      +\frac{1}{C\left(\frac{T-t}{T}\right)^2
                +\frac{t(T-t)}{T}}
    =\frac{1}{T-t}\left(-1+\frac{T^2}{C(T-t)+tT}\right)\\
    &=\frac{-C(T-t)+T(T-t)}{(T-t)(C(T-t)+tT)}
    =\frac{T-C}{(T-C)t+CT},
    \qquad t\in[0,T).
\end{align*}
Hence
 $$
  q_C(t)=
     \begin{cases}
          0 & \text{\quad if \ $C=T,$}\\[1mm]
          \frac{1}{t+\frac{CT}{T-C}}
             & \text{\quad if \ $C\not=T,$}
     \end{cases}
     \qquad \text{for all}\quad
     t\in[0,T).
 $$
In case of \ $C=T$ \ (and \ $\alpha(t)=1$, $t\in[0,T)$), \ the process \ $Z$ \ in Theorem \ref{THEOREM5}
 is a standard Wiener process, which corresponds to the case \ $C=+\infty$ \ in Theorem \ref{THEOREM6}.
Moreover, since the range of the function
 \ $(0,+\infty)\setminus\{T\}\ni C\mapsto\frac{CT}{T-C}$ \ is
 \ $(-\infty,-T)\cup(0,+\infty),$ \ we have the family of the given
 functions \ $q_C,$ $C>0,$ \ can be parametrized also in the
 form
 $$
   q_{\widetilde C}(t)=\frac{1}{t+\widetilde C},
    \qquad t\in[0,T),
 $$
where \ $\widetilde C\in(\RR\setminus[-T,0])\cup\{+\infty\}.$ \
Moreover, since \ $\lim_{t\uparrow T}q_{\widetilde C}(t)=(T+\widetilde C)^{-1}$, \ the assumptions of part
 (ii) of Theorem \ref{THEOREM5} are also satisfied.
Hence we get Theorem \ref{THEOREM6}.
\proofend

\begin{Rem}
Note that for a fixed \ $T>0$ \ and \ $C\in\RR\setminus[-T,0],$ \
there are many continuous functions \ $q:[0,\infty)\to\RR$ \ for which
 \ $q(t)=\frac{1}{t+C},$ $t\in[0,T).$ \
Note also that for \ $C>0$ $(C\in\RR)$ \ or \ $C=+\infty$ \ the
 continuously differentiable function \ $q(t)=\frac1{t+C},\,t\geq0,$ \ does not depend
 on \ $T$ \ and hence in this case for all \ $T>0$ \ the law of
 the Ornstein-Uhlenbeck type bridge from $0$ to $0$ over \ $[0,T]$ \ derived
 from \ $Z$ \ (with this function \ $q$ \ and \ $\sigma(t)=1$, $t\in[0,T)$, \ or \ $\sigma(t)=-1$, $t\in[0,T)$)
 coincides with the law of the usual Wiener bridge from $0$ to $0$ over \ $[0,T].$
 \proofend
\end{Rem}

\par\bigskip\noindent
{\bf Acknowledgment.} We would like to thank Endre Igl\'oi
 for some very helpful remarks concerning Lemma \ref{LEMMA3}.


\begin{thebibliography}{99}

\bibitem{BarIgl} {\sc M. Barczy} and {\sc E. Igl\'oi},
Karhunen-Lo\`{e}ve expansions of alpha-Wiener bridges.
{\sl Central European Journal of Mathematics} {\bf 9}(1)
 (2011), 65-84.

\bibitem{BarBec} {\sc M. Barczy} and {\sc P. Kern},
 Representations of multidimensional linear process bridges.
{\em Arxiv}, URL: {\texttt http://arxiv.org/abs/1011.0067}

\bibitem{BarPap2} {\sc M. Barczy} and {\sc G. Pap},
 Alpha-Wiener bridges: singularity of induced measures
 and sample path properties.
{\sl Stochastic Analysis and Applications} {\bf 28}(3)
   (2010), 447-466.

\bibitem{Bau} {\sc H. Bauer},
   {\sl Probability Theory}. Walter de Gruyter, 1996.

\bibitem{BenLee} {\sc I. Benjamini} and {\sc S. Lee}, Conditioned diffusions
  which are Brownian bridges.
 {\sl Journal of Theoretical Probability} {\bf 10}(3)
   (1997), 733-736.

\bibitem{Bor} {\sc A. N. Borodin},
  Diffusion Processes with Identical Bridges.
  {\sl Journal of Mathematical Sciences} {\bf 127}(1) (2005), 1687--1695.

\bibitem{BreSch} {\sc M. J. Brennan} and {\sc E. S. Schwartz},
    Arbitrage in stock index futures. {\sl The Journal of Business} {\bf 63}(1)
   (1990), S7-S31.

\bibitem{CsorRev} {\sc M. Cs\"org\H o} and {\sc P. R\'ev\'esz},
{\sl Strong Approximations in Probability and Statistics.}
Academic Press, New York, 1981.

\bibitem{DzhSpr} {\sc K. Dzhaparidze} and {\sc P. Spreij},
 The strong law of large numbers for martingales with deterministic quadratic variation.
 {\sl Stochastics and Stochastics Reports}  {\bf 42}(1)  (1993), 53--65.

\bibitem{Fit} {\sc P. J. Fitzsimmons}, Markov processes with identical bridges.
    {\sl Electronic Journal of Probability} {\bf 3} (1998), Paper no.
    12., 1-12.

\bibitem{JacShi} {\sc J. Jacod} and  {\sc A. N. Shiryaev},
   Limit Theorems for Stochastic Processes, 2nd edition. Springer-Verlag, Berlin, 2003.


\bibitem{KarShr} {\sc I. Karatzas} and {\sc S. E. Shreve},
    {\sl Brownian Motion and Stochastic Calculus}, 2nd edition.
    Springer-Verlag, Berlin, Heidelberg, 1991.

\bibitem{Kov} {\sc V. A. Koval'}, On the strong law of large numbers for multivariate martingales
 with continuous time. {\sl Ukrainian Mathematical Journal} {\bf 53}(9) (2001), 1554--1560.

\bibitem{Man} {\sc R. Mansuy}, On a one-parameter generalization of the Brownian bridge
   and associated quadratic functionals.
   {\sl Journal of Theoretical Probability} {\bf 17}(4) (2004), 1021--1029.

\bibitem{Oks} {\sc B. {\O}ksendal},
 {\sl Stochastic Differential Equations, 6th edition.}
 Springer, Berlin, Heidelberg, New York, 2003.

\bibitem{Rei} {\sc W. T. Reid},
 {\sl Riccati Differential Equations}.
 Academic Press, New York and London, 1972.

\bibitem{RevYor} {\sc D. Revuz} and {\sc M. Yor},
 {\sl Continuous martingales and Brownian motion}, 3rd edition,
 corrected 2nd printing. Springer-Verlag, Berlin, 2001.

\bibitem{Rog}
{\sc L. C. G. Rogers},
 Smooth transition densities for one-dimensional diffusions.
 {\sl The Bulletin of the London Mathematical Society } {\bf 17}(2) (1985), 157--161.

\bibitem{Shi} {\sc A. N. Shiryaev}, {\sl Probability,} 2nd edition. Springer, 1989.

\bibitem{SonTreWil}
{\sc D. Sondermann, M. Trede} and {\sc B. Wilfling},
Estimating the degree of interventionist policies in the run-up to EMU.
{\sl Applied Economics} {\bf 43}(2) (2011), 207--218.

\bibitem{TreWil}
{\sc M. Trede} and {\sc B. Wilfling},
Estimating exchange rate dynamics with diffusion processes: an application to Greek EMU data.
 {\sl Empirical Economics} {\bf33}(1) (2007), 23--39.

\bibitem{Vra} {\sc I. I. Vrabie}, {\sl Differential Equations,
 An Introduction to Basic Concepts, Results and Applications.}
 World Scientific Publishing Co. Pte. Ltd, 2004.

\bibitem{Wal} {\sc W. Walter}, {\sl Gew\"ohnliche Differentialgleichungen.}  Springer, 2000.

\end{thebibliography}
\end{document}